\documentclass[11pt, a4paper]{article}

\usepackage[T1]{fontenc}
\usepackage[utf8]{inputenc}
\usepackage{xcolor}
\usepackage{graphicx}
\usepackage{times}
\usepackage{amsmath}
\usepackage{amssymb}
\usepackage{subcaption}

\usepackage{url}
\usepackage{multirow}
\usepackage{cite}
\usepackage{algorithm}

\usepackage{algorithmic}
\usepackage{microtype}
\newcommand{\eqdef}{\triangleq}
\newcommand{\qedsymbol}{\hfill$\square$}
\newcommand{\dvec}{d}

\newcommand{\diag}{\operatorname{diag}}

\newtheorem{theorem}{Theorem}[section]
\newtheorem{lemma}[theorem]{Lemma}

\newtheorem{proposition}[theorem]{Proposition}

\newtheorem{assumption}[theorem]{Assumption}

\newenvironment{proof}{\emph{Proof.}}{\hfill\qedsymbol}

\newcommand{\rr}{\mathbb{R}}

\begin{document}

\title{A Parametric Active-Set Algorithm for Solving Monotone
Generalized Nash Equilibrium Problems}

\author{Alberto Bemporad and Tatiana Tatarenko 
{\renewcommand{\thefootnote}{}\thanks{A. Bemporad is with the IMT School for Advanced Studies, Piazza San Francesco 19, Lucca, Italy. Email: \texttt{\scriptsize alberto.bemporad@imtlucca.it}.
T. Tatarenko is with the Department of Control Theory and Intelligent Systems, TU Darmstadt, Germany. E-mails: \texttt{\scriptsize tatiana.tatarenko@tu-darmstadt.de}.\\   
This work was funded by the European Union (ERC Advanced Research Grant COMPACT, No. 101141351). Views and opinions expressed are however those of the authors only and do not necessarily reflect those of the European Union or the European Research Council. Neither the European Union nor the granting authority can be held responsible for them. The work was also funded by the Deutsche Forschungsgemeinschaft (DFG, German Research Foundation). Project number: 528033031.
}}}

\maketitle
\thispagestyle{empty}

\begin{abstract}
We present a fast active-set algorithm for the solution of monotone linear-quadratic generalized Nash equilibrium problems (GNEPs) with shared affine inequality and equality constraints. The algorithm, which we call gPAQ, is an adaptation of the parametric active-set method proposed by Best~(1996) and extended in qpOASES by Ferreau et al.~(2014) to solve convex quadratic programs. Under strong monotonicity, we show that the same active-set strategy applies. We handle the case of merely monotone problems by using gPAQ
as the inner solver within proximal-point iterations. 
To solve strongly monotone convex nonlinear GNEPs, we propose a Newton-type iterative scheme, which we call gPAN, which iteratively linearizes the nonlinear problem and solves the resulting LQ-GNEP subproblems using gPAQ.
We show in numerical examples that the proposed algorithms are significantly faster than existing methods.
\end{abstract}

\section{Introduction}
Generalized Nash equilibrium problems (GNEPs) provide a natural framework for multi-agent decision-making in which strategic agents have objectives and feasible strategies coupled with those of the other players. Such problems arise in different engineering applications, including power and communication networks, and autonomous vehicles. In particular, in game-theoretic model predictive control (GT-MPC) approaches, a GNEP must be solved at every sampling instant~\cite{HBLD22}. Therefore, the practical impact of game-theoretic control approaches hinges on the availability of solvers that compute a generalized Nash equilibrium (GNE) reliably and fast, a need that has motivated the development of dedicated software tools~\cite{LSM22,Bem26,ABB26}.

An important and widely studied subclass of GNEPs is obtained when each player minimizes a convex quadratic cost and the shared coupling constraints are affine, referred to as linear-quadratic GNEPs (LQ-GNEPs).
The main applications include linear GT-MPC formulations~\cite{HB26,BG26,DFB26} and LQ-GNEPs learned from data~\cite{BT26c}.
Under the standard assumption that the game's stacked pseudogradient is monotone, a distinguished equilibrium, the variational GNE (v-GNE), can be characterized as the solution of a monotone affine variational inequality (VI) over the shared feasible set~\cite{FP03}. 

The v-GNE characterization underlies most existing solution methods, typically based on forward-backward, proximal-point, or extragradient-type splitting schemes~\cite{BYGP22,Kor76}, Douglas-Rach\-ford splitting methods like DR-DAQP~\cite{ABB26}, the Operator Extrapolation method~\cite{Lan2022} on strongly monotone instances, and Korpelevich's extragradient method~\cite{Kor76} on merely monotone instances.
While theoretically appealing, many of these iterative methods only attain fast, linearly convergent behavior when the pseudogradient is \emph{strongly} monotone~\cite{TatKam25}.
For convex nonlinear monotone GNEPs, a general-purpose nonlinear least-squares solver can be employed to minimize the game's KKT residual, as implemented for instance in NashOpt~\cite{Bem26}. 

A different line of work exploits the linearity of the joint Karush--Kuhn--Tucker (KKT) conditions through pivoting, in the spirit of Lemke's classical algorithm for linear complementarity problems, either specialized to LQ-GNEPs with shared constraints~\cite{SPS13} or applied more generally to affine VIs~\cite{CF96,KHF18}, or by recasting the LQ-GNEP as a QP that can be solved by standard QP solvers~\cite{BT26b}. Such methods compute an exact GNE, when one exists, after a finite number of pivots, without resorting to an iterative fixed-point argument. 

In this paper, we also consider pivoting algorithms by investigating whether Best's parametric active-set method~\cite{Bes96} can be extended to solve LQ-GNEPs. The method, extended and popularized by the qpOASES solver~\cite{FKPBD14}, follows a related, finitely terminating pivoting philosophy, but through a different mechanism: it tracks the KKT solution of a convex QP along a homotopy that interpolates between a trivial problem and the QP of interest, updating the optimal active set only at a finite number of breakpoints. Since its derivation only relies on the positive definiteness of the (symmetric) QP Hessian, it is natural to ask whether the same active-set strategy extends to the generally \emph{asymmetric} pseudogradient matrix of an LQ-GNEP.

We answer this question affirmatively. The main contribution of this paper is an algorithm, which we call gPAQ (``\textbf{g}eneralized Nash equilibrium-seeking \textbf{P}arametric \textbf{A}ctive-set method for linear-\textbf{Q}uadratic games''), that extends Best's parametric active-set homotopy to LQ-GNEPs with shared affine inequalities, and, after a simple elimination step, equality constraints if present. We show that, in spite of the possible asymmetry of the stacked pseudogradient matrix, the same active-set path-following strategy and its finite-termination guarantee carry over intact whenever the game is strongly monotone (Sections~\ref{sec:problem-formulation} and~\ref{sec:active-set-path-following})
and discuss how to support warm starting (Section~\ref{sec:warm-start}) and equality constraints (Section~\ref{sec:equality}).
We then handle \textit{merely monotone} LQ-GNEPs, for which strong-monotonicity-based guarantees no longer apply, by using gPAQ as the inner solver of a proximal-point outer loop, warm-started at each iteration from the previously optimal active set (Section~\ref{sec:merely-monotone}). Finally, adopting a Newton-type outer iteration in the spirit of~\cite{TFI93}, we propose a method, called gPAN, to solve convex \textit{nonlinear} strongly monotone GNEPs by repeatedly linearizing the pseudogradient and solving the resulting LQ-GNEP subproblem with gPAQ (Section~\ref{sec:nl-monotone}). Section~\ref{sec:examples} shows the superior performance of gPAQ and its extensions compared to several existing state-of-the-art solvers and Section~\ref{sec:conclusion} concludes the paper. A Python/C implementation of the gPAQ and gPAN algorithms is available in the \texttt{nashopt} library~\cite{Bem26},
the code to reproduce the reported tests is available at \texttt{github.com/bemporad/gpaq\_examples}.

\subsection{Notation}
We denote by $\rr^n$ the $n$-dimensional Euclidean space. Given two vectors $a\in\rr^n$, $b\in\rr^m$, we denote by $(a,b)$ their concatenation $[a'\ b']'\in \rr^{n+m}$.
Given a finite set $J=\{j_1,\ldots,j_K\}\subseteq\{1,\ldots,m\}$, $|J|=K\leq m$ denotes its cardinality and $\bar J=\{1,\ldots,m\}\setminus J$; for a matrix $A\in\rr^{m\times n}$ and a vector $b\in\rr^m$, $A_J$ and $b_J$ denote the submatrix of the rows of $A$ and the subvector of the entries of $b$ indexed by $J$, respectively, and $a_j'$ (also written $A_j$) denotes the $j$th row of $A$. For $x=(x_1,\ldots,x_N)$, we denote $x_{-i}\eqdef(x_1,\ldots,x_{i-1},x_{i+1},\ldots,x_N)$. We denote by $I$ the identity matrix. For a square matrix $M$, $M\succ0$ ($M\succeq0$) means $x'Mx>0$ for all $x\neq0$ ($x'Mx\geq0$ for all $x$), and $M\succ N$ ($M\succeq N$) means $M-N\succ0$ ($M-N\succeq0$); if $M$ is symmetric, $\lambda_{\rm min}(M)$ denotes its smallest eigenvalue.

\section{Problem formulation and preliminaries}
\label{sec:problem-formulation}
Consider a game with $N$ agents deciding a vector $x_i$ of decision variables,
$x_i\in\rr^{n_i}$, and let $x=(x_1,\ldots,x_N)\in\rr^n$, $n=\sum_i n_i$. 
Player $i$ solves the quadratic program
\begin{equation}
\label{eq:player}
\begin{aligned}
 \min_{x_i}& \frac12 x_i'C_{ii}x_i +x_i'\sum_{j\ne i}C_{ij}x_j+c_i'x_i\\
\textrm{s.t. }& Ax\le b
\end{aligned}
\end{equation}
where the affine inequality constraints are shared by all players, $A\in\rr^{m\times n}$, $b\in\rr^m$. 
We assume that
$C_{ii}=C_{ii}'\succ 0$, $\forall i=1 \ldots N$, and let $K\eqdef \{x\in\rr^n:\ Ax\leq b\}$.
The stacked pseudogradient $F(x)\eqdef Cx+c$ associated with~\eqref{eq:player} is used
to characterize the {\it variational GNE} (v-GNE) $x^\star\in K$ of the game
and the common multiplier vector $\lambda^\star\ge0$, $\lambda^\star\in\rr^m$, 
through the following joint KKT conditions~\cite{FK07}:
\begin{equation}
\begin{aligned}
 &Cx^\star+c+A'\lambda^\star=0,&&\lambda^{\star\prime}(Ax^\star-b)=0\\
 &Ax^\star\le b,&& \lambda^\star\geq 0.
\label{eq:kkt}
\end{aligned}
\end{equation}

\begin{assumption}[Strong monotonicity]\label{ass:strong}
There exists a scalar $\mu>0$ such that
\begin{equation}
 C_s\eqdef \frac12(C+C')\succeq \mu I.
\label{eq:strong-monotonicity}
\end{equation}
\end{assumption}
Note that Assumption~\ref{ass:strong} implies $x'Cx=x'C_sx\geq\mu\|x\|^2>0$, $\forall x\neq0$, 
which also implies that $C$ is invertible.
Recall that $x^\star$ is a v-GNE of~\eqref{eq:player} if and only if it solves the variational inequality VI$(K,F)$, i.e., $x^\star\in K$ and $F(x^\star)'(y-x^\star)\ge0$ for all $y\in K$~\cite{FK07}, which, since $K$ is polyhedral, is in turn equivalent to the existence of $\lambda^\star$ such that $(x^\star,\lambda^\star)$ satisfies~\eqref{eq:kkt}. Under Assumption~\ref{ass:strong}, the pseudogradient $F$ is strongly monotone, as $(F(x)-F(y))'(x-y)=(x-y)'C_s(x-y)\ge\mu\|x-y\|^2$ for all $x,y\in\rr^n$. Hence, VI$(K,F)$ has at most one solution, and exactly one whenever $K\neq\varnothing$~\cite[Thm.~2.3.3]{FP03}, i.e., the game~\eqref{eq:player} admits a unique v-GNE.

\subsection{Parametric homotopy}
Similarly to~\cite{FKPBD14}, we consider the GNEP defined in~\eqref{eq:player} as a special instance of a parametric family GNEP$(t)$, where $t\in[0,1]$ is a homotopy parameter, GNEP$(1)$ coincides with the original problem~\eqref{eq:player}, GNEP$(0)$ is a simple problem with a known solution, and the solution of GNEP$(t)$ is tracked as $t$ increases from $0$ to $1$,
adapting the approach of~\cite{Bes96}.

For an arbitrary vector $\dvec\in\rr^m$ with positive entries, define the following homotopy:
\begin{equation}\label{eq:homotopy}
\textrm{GNEP}(t):\left\{\begin{array}{rl}
\displaystyle{\min_{x_i}}& \displaystyle{\frac12 x_i'C_{ii}x_i +x_i'\sum_{j\ne i}C_{ij}x_j+c_i(t)'x_i}\\
\text{s.t. } & Ax\le b(t)
\end{array}\right.
\end{equation}
where 
\[
    b(t)\eqdef \dvec+tp \qquad p\eqdef (b-\dvec) \qquad c(t)\eqdef tc.
\]
Denote by $K(t)\eqdef \{x:\ Ax\leq b(t)\}$ its associated feasible set,
and by $F_t(x)\eqdef Cx+tc$ the pseudogradient associated with GNEP$(t)$. We denote by $x(t)$ the v-GNE of GNEP$(t)$ and by $\lambda(t)$ the associated common multiplier, with $x^\star=x(1)$ being the v-GNE of~\eqref{eq:player} we want to compute. 

\subsection{Affine solution for a fixed active set}
The functional form of $(x(t),\lambda(t))$ is determined by the corresponding optimal active set $J(t)$, as described in the following.

\begin{proposition}
\label{prop:bordered}
Let Assumption~\ref{ass:strong} hold, and consider a generic index set $J\subseteq\{1 \dots m\}$ such that the corresponding rows of $A$
form a full row-rank matrix $A_J$, and define
\begin{equation}\label{eq:bordered}
    H_J\eqdef \begin{bmatrix}C&A_J'\\A_J&0\end{bmatrix}.
\end{equation}
Then $H_J$ is nonsingular.
\end{proposition}
\begin{proof}
The proof is straightforward. For given vectors $s$ and $\xi$ such that $H_J[s'\ \xi']'=0$, we get $Cs+A_J'\xi=0$, $A_Js=0$, and therefore, premultiplying by $s'$, we get $0=s'Cs+(A_Js)'\xi = s' C_s s\geq\mu\|s\|^2$, so necessarily $s=0$. This in turn gives $A_J'\xi=0$, and full row-rank of $A_J$ gives $\xi=0$.
\end{proof}

Assuming that $A_J$ has full row rank, by Proposition~\ref{prop:bordered}, we can solve the following two linear systems
\begin{subequations}
\begin{equation}\label{eq:two-systems}
 H_J\begin{bmatrix}h_{1}\\g_{1J}\end{bmatrix}
 =\begin{bmatrix}0\\\dvec_J\end{bmatrix}
 \qquad
 H_J\begin{bmatrix}h_{2}\\g_{2J}\end{bmatrix}
 =\begin{bmatrix}-c\\p_J\end{bmatrix}
\end{equation}
and define the affine path
\begin{equation}
\label{eq:affine-path-x-lambda}
 x(t)=h_{1}+h_{2}t,\quad
 \lambda_{J}(t)=g_{1J}+g_{2J}t,\quad \lambda_{\bar J}(t)=0.
\end{equation}
\label{eq:affine-path}%
\end{subequations}

\begin{proposition}
\label{prop:affine}
For every $t$ such that $A_{J(t)} x(t)=b_{J(t)}(t)$, $A_{\bar J(t)} x(t)\leq b_{\bar J(t)}(t)$, $\lambda_{J(t)}(t)\geq 0$, and $\lambda_{\bar J(t)}(t)=0$, the affine pair $(x(t),\lambda(t))$ defined in \eqref{eq:affine-path} solves GNEP$(t)$, and $x(t)$ is the unique v-GNE of GNEP$(t)$.
\end{proposition}
\begin{proof}
By adding the first system in \eqref{eq:two-systems} to $t$ times the second, we obtain that $x(t),\lambda(t)$ satisfy
\begin{equation}
 Cx(t)+tc+A_J'\lambda_{J}(t)=0
 \qquad A_J x(t)=b_J(t).
\label{eq:stationarity-t}
\end{equation}
The remaining KKT conditions of primal and dual feasibility are satisfied for all $t$ such that the
inequalities $A_{\bar J} x(t)\leq b_{\bar J}(t)$ and $\lambda_J(t)\geq 0$ hold, so that $(x(t),\lambda(t))$
solves GNEP$(t)$ for such $t$. Uniqueness of the v-GNE solution $x(t)$ follows from $C_s\succ 0$ and 
$K(t)\neq \emptyset$~\cite[Thm.~2.3.3]{FP03}.
\end{proof}

The following simple result will be used as in~\cite{Bes96} to advance the parameter $t$ along the path while preserving strict feasibility of the affine solution $(x(t),\lambda(t))$.
\begin{proposition}
\label{prop:persistence}
Let $v(t)=v_0+v_1(t-\tau)\in\rr^k$, $v_0,v_1\in\rr^k$, $\tau\in\rr$, be an affine function of $t\in\rr$ and $v_0>0$ componentwise. Then $v(t)>0$ for every $t\in[\tau,\tau+\epsilon)$ where
\[
    \epsilon=\min\{-v_{0i}/v_{1i}:\ v_{1i}<0\}
\]
with $\epsilon=+\infty$ in the case $v_{1i}\geq 0$ for all $i=1,\ldots,k$.
\end{proposition}
\begin{proof}
Fix a component $i$. If $v_{1i}\geq0$ then $v_i(t)=v_{0i}+v_{1i}(t-\tau)\geq v_{0i}>0$ for every $t\geq\tau$. If $v_{1i}<0$ then $v_i(t)=v_{0i}+v_{1i}(t-\tau)>0$ for every $t\in[\tau,\tau+\epsilon_i)$, where $\epsilon_i\eqdef -v_{0i}/v_{1i}>0$. Taking $\epsilon\eqdef \min_i\epsilon_i$ over the indices with $v_{1i}<0$ (or $\epsilon\eqdef +\infty$ if none) gives $v(t)>0$ for every $t\in[\tau,\tau+\epsilon)$.
\end{proof}

\section{Active-set path-following}
\label{sec:active-set-path-following}
Our goal is to extend the analysis of~\cite{Bes96,FKPBD14} to the GNEP case, for which $C$ is not symmetric in general, and show how to track the optimal active set $J(t)$ along the path from $t=0$ to $t=1$. We will show that the path is piecewise-affine in that the active set changes only at a finite number of breakpoints $t_1,\ldots,t_K$, $t_K\geq 1$ (with possibly $t_K=+\infty$), where either an inactive constraint becomes tight ($A_jx(t)=b_j(t)$ for some $j\in\bar J$) or an active multiplier 
$\lambda_j(t)$ becomes zero for some $j\in J(t)$. 

For $t=0$, $x(0)=0$ and $\lambda(0)=0$ solve GNEP$(0)$, with the associated optimal active set $J(0)=\varnothing$
since $0=Ax(0)<d=b(0)$. The next proposition shows that the path also 
remains strictly feasible in a sufficiently small interval $(0,\epsilon_0)$.
\begin{proposition}
Assume that $C$ is invertible and $d>0$. Then, there exists a scalar $\epsilon_0>0$ (possibly $+\infty$) such that, for all $t\in[0,\epsilon_0)$, $x(t)=-tC^{-1}c$ and $\lambda(t)=0$ solve GNEP$(t)$, and $x(t)$ is strictly feasible.
\label{prop:start}
\end{proposition}
\begin{proof}
The KKT stationarity condition $Cx(t)+tc+A'\lambda(t)=0$ is satisfied by $x(t)=-tC^{-1}c$, $\lambda(t)=0$, where
the latter implies that the complementarity condition and dual feasibility $\lambda(t)\geq 0$ are also satisfied. The only remaining condition to check is strict primal feasibility $Ax(t)< b(t)$, which is equivalent to $v(t)\eqdef d+(AC^{-1}c+p)t>0$. Since $d>0$, by Proposition~\ref{prop:persistence}, there exists $\epsilon_0>0$ (or $\epsilon_0=+\infty$) such that $v(t)>0$ for all $t\in(0,\epsilon_0)$, so that $x(t)$ 
remains strictly feasible for all $t\in(0,\epsilon_0)$.
\end{proof}

Note that, as observed earlier, Assumption~\ref{ass:strong} implies that $C$ is invertible.

Consider now a generic active set $J_j\subseteq\{1,\ldots,m\}$ such that
$A_{J_j}$ has full row rank and such that there exist $t_j\ge0$ and $\epsilon_j>0$ for which the pair
$(x(t),\lambda(t))$ obtained by solving~\eqref{eq:affine-path} for $J=J_j$ satisfies
\[
\begin{aligned}
 A_{J_j}x(t)&=b_{J_j}(t),\quad \lambda_{\bar J_j}(t)=0,&& t\in[t_j,t_j+\epsilon_j)\\
 \lambda_{J_j}(t)&>0,\quad A_{\bar J_j}x(t)<b_{\bar J_j}(t),&& t\in(t_j,t_j+\epsilon_j).
\end{aligned}
\]
Let $x(t)=h_{1j}+h_{2j}t$ and $\lambda_{J_j}(t)=g_{1j}+g_{2j}t$ denote
the corresponding solutions obtained from~\eqref{eq:two-systems} and define
\begin{subequations}
\begin{equation}\label{eq:primal-event}
 \widehat t_{j+1}\eqdef 
 \min_{\substack{i\notin J_j\\a_i'h_{2j}>p_i\\r_i>t_j}}
 r_i,
 \qquad
 r_i\eqdef \frac{\dvec_i-a_i'h_{1j}}{a_i'h_{2j}-p_i}
\end{equation}
corresponding to the first inactive constraint that reaches equality, and
\begin{equation}\label{eq:dual-event}
 \widetilde t_{j+1}\eqdef 
 \min_{\substack{i\in J_j\\(g_{2j})_i<0\\q_i>t_j}}
 q_i,
 \qquad
 q_i\eqdef -\frac{(g_{1j})_i}{(g_{2j})_i}
\end{equation}
\label{eq:events}%
\end{subequations}
corresponding to the first active multiplier that reaches zero, 
where we have
assumed that the minimum of an empty set is $+\infty$. 
Let
$t_{j+1}\eqdef \min\{\widehat t_{j+1},\widetilde t_{j+1}\}$. As in~\cite{Bes96},
we distinguish different cases,
depending on the values of $\widehat t_{j+1}$ and $\widetilde t_{j+1}$.

\paragraph*{Case 0: $\widehat t_{j+1}=\widetilde t_{j+1}=\infty$.}~
The parametric solution $(x(t),\lambda(t))$ satisfies all the joint KKT conditions and, therefore, solves GNEP$(t)$ for all $t\ge t_j$, in particular for $t=1$.

\paragraph*{Case 1: $\widehat t_{j+1}<\widetilde t_{j+1}$.}~
A constraint $h\in\bar J_j=\{1,\ldots,m\}\setminus J_j$ becomes active at $t_{j+1}=\widehat t_{j+1}$, so we set $J_{j+1}\eqdef J_j\cup\{h\}$. Let us assume for now that all the remaining constraints in $\bar J_j$ 
are inactive at $t_{j+1}$ (i.e., the arg min in~\eqref{eq:primal-event} is unique) and that $A_{J_{j+1}}$ has full row rank (i.e., the new active row $a_h$ is linearly independent of the rows of $A_{J_j}$).

\paragraph*{Case 2: $\widehat t_{j+1}>\widetilde t_{j+1}$.}~
A multiplier $\ell\in J_j$ becomes zero at $t_{j+1}=\widetilde t_{j+1}$, so we set $J_{j+1}\eqdef J_j\setminus\{\ell\}$. Let us assume for now that all the remaining multipliers in $J_j$ are positive at $t_{j+1}$ (i.e., the arg min in~\eqref{eq:dual-event} is unique).

As in~\cite[p. 63]{Bes96}, besides the uniqueness of the arg min in~\eqref{eq:primal-event} and~\eqref{eq:dual-event}, we have also assumed that $\widehat t_{j+1}=\widetilde t_{j+1}<\infty$ does not occur. We will discuss these (rare) events in Section~\ref{sec:relaxing-uniqueness}.

We now show that, in both Case 1 and Case 2, the affine solution $x(t)$ associated with $J_{j+1}$ is strictly feasible and the multipliers $\lambda_{J_{j+1}}(t)$ are strictly positive immediately to the right of $t_{j+1}$, so that the parameter $t$ can be further increased.

\begin{lemma}%
\label{lem:pivot-sign}
Let Assumption~\ref{ass:strong} hold, let $J\subseteq\{1,\ldots,m\}$ be such that $A_J$ has full row rank, and let $a\in\rr^n$ be linearly independent of the rows of $A_J$. Let $z=[z_x'\ z_\lambda']'$ solve
\begin{equation}
H_J\begin{bmatrix}z_x\\z_\lambda\end{bmatrix}=
\begin{bmatrix}C&A_J'\\A_J&0\end{bmatrix}\begin{bmatrix}z_x\\z_\lambda\end{bmatrix}=\begin{bmatrix}a\\0\end{bmatrix}.
\label{eq:z}
\end{equation}
Then $a'z_x>0$.
\end{lemma}
\begin{proof}
If $z_x=0$ then from~\eqref{eq:z} we get $A_J'z_\lambda=a$, i.e., $a$ is a linear combination of the rows of $A_J$, 
contradicting the hypothesis, so $z_x\neq 0$.
Premultiplying $Cz_x+A_J'z_\lambda=a$ by $z_x'$ and using $A_Jz_x=0$ we get $a'z_x=z_x'Cz_x=z_x'C_sz_x\geq\mu\|z_x\|^2>0$
by Assumption~\ref{ass:strong}.
\end{proof}

\begin{lemma}
In both Case 1 and Case 2, there exists a scalar $\epsilon_{j+1}>0$ such that the affine solution $(\bar x(t),\bar\lambda(t))$ associated with $J_{j+1}$ given by~\eqref{eq:affine-path} satisfies
\begin{subequations}
\begin{equation}
    A_{J_{j+1}}\bar x(t)=b_{J_{j+1}}(t),\quad \bar\lambda_{\bar J_{j+1}}(t)=0
    \label{eq:active-equalities}
\end{equation}
for all $t\in[t_{j+1},t_{j+1}+\epsilon_{j+1})$, and 
\begin{equation}
     \bar\lambda_{J_{j+1}}(t)>0,\quad A_{\bar J_{j+1}}\bar x(t)<b_{\bar J_{j+1}}(t)
    \label{eq:inactive-inequalities}
\end{equation}%
\label{eq:active-inactive-equalities}%
\end{subequations}
for all $t\in(t_{j+1},t_{j+1}+\epsilon_{j+1})$.
\label{lem:strict-after-pivot}
\end{lemma}
\begin{proof}
The proof is similar to that reported in~\cite{Bes96}, see Appendix~\ref{app:strict-after-pivot}.
\end{proof}

\subsection{Relaxing the independence assumption}
\label{sec:relaxing-independence}
In Case 1 and in Lemma~\ref{lem:strict-after-pivot}, we assumed the linear independence of $a_h$ from the rows of $A_{J_j}$. We now relax this assumption.

\begin{lemma}
In the setting of Case 1, suppose $a_h=A_{J_j}'\xi$ for some $\xi\in\rr^{|J_j|}$. Then:

\noindent ($i$) If $\xi_i>0$ for some $i\in J_j$, let
\[
 k\eqdef \arg\min_{i\in J_j:\ \xi_i>0}\frac{\lambda_{i}(t_{j+1})}{\xi_i},\quad \theta^\star=\frac{\lambda_{k}(t_{j+1})}{\xi_k}
\]
and set $J_{j+1}\eqdef (J_j\setminus\{k\})\cup\{h\}$. Then $A_{J_{j+1}}$ has full row rank, and
\begin{equation}
\begin{aligned}
&\bar x(t_{j+1})\eqdef x(t_{j+1}),\quad
 \bar\lambda_h(t_{j+1})\eqdef \theta^\star\\
& \bar\lambda_i(t_{j+1})\eqdef \lambda_{J_j,i}(t_{j+1})-\theta^\star\xi_i,\ i\in J_j\setminus\{k\}
\end{aligned}
\label{eq:switch}
\end{equation}
is a KKT pair for the active set $J_{j+1}$ at $t=t_{j+1}$, with $\bar\lambda_{J_{j+1}}(t_{j+1})\geq 0$.

\noindent ($ii$) If $\xi_i\leq 0$ for every $i\in J_j$, then GNEP$(t)$ has no feasible solution for any $t>t_{j+1}$ and in particular, if $t_{j+1}<1$, the given GNEP~\eqref{eq:player} has no feasible solution due to the emptiness of the set $\{x:\ Ax\leq b\}$.
\label{lem:dependent-exchange}
\end{lemma}

\begin{proof}
The proof is similar to that reported in~\cite{Bes96}, see Appendix~\ref{app:dependent-exchange}.
\end{proof}

Case ($ii$) of Lemma~\ref{lem:dependent-exchange} is the exact analogue of the infeasibility branch of~\cite[Case 1.2.2]{Bes96}.

\begin{lemma}
In the setting of case ($i$) of Lemma~\ref{lem:dependent-exchange}, suppose in addition that the minimizer $k$ in~\eqref{eq:switch} is unique. Then there exists $\epsilon_{j+1}>0$ such that the affine pair $(\bar x(t),\bar\lambda(t))$ of~\eqref{eq:switch} satisfies
$\bar\lambda_{J_{j+1}}(t)>0$, $A_{\bar J_{j+1}}\bar x(t)<b_{\bar J_{j+1}}(t)$, $\forall t\in(t_{j+1},t_{j+1}+\epsilon_{j+1})$.
\label{lem:strict-after-dependent-exchange}
\end{lemma}
\begin{proof}
The proof is similar to that reported in~\cite{Bes96}, see Appendix~\ref{app:strict-after-dependent-exchange}.
\end{proof}

The following lemma shows that the same active set $J$ cannot be visited twice along the pivot sequence.
\begin{lemma}
\label{lem:no-repeat}
Any active set $J_j\subseteq\{1,\ldots,m\}$ with $A_{J_j}$ full row rank which is active
on $[t_j,t_{j+1}]$ cannot be active again on $[t_k,t_{k+1}]$ along the homotopy path if $t_k>t_{j+1}$.
\end{lemma}

\begin{proof}
The proof is similar to that reported in~\cite{Bes96}, see Appendix~\ref{app:no-repeat}.
\end{proof}

\subsection{Main result}
We are ready to state our main result, which is the analogue of~\cite[Thm.~4.1]{Bes96} for the GNEP~\eqref{eq:player},
when running the homotopy~\eqref{eq:homotopy} with the parametric active-set strategy described above until either some $t_{j+1}\ge1$ is reached or infeasibility is detected.
\begin{theorem}%
Let Assumption~\ref{ass:strong} hold, and let $\dvec>0$ be such that, at every breakpoint $t_{j+1}$,
the following nondegeneracy conditions hold:
($a$) whichever of $\widehat t_{j+1}$, $\widetilde t_{j+1}$ in~\eqref{eq:events} is finite has a unique minimizer; ($b$) $\widehat t_{j+1}\neq\widetilde t_{j+1}$ whenever both are finite; and ($c$) whenever case~($i$) of Lemma~\ref{lem:dependent-exchange} is invoked, the minimizer $k$ in~\eqref{eq:switch} is likewise unique.
Starting from $J_0=\varnothing$, $t_0=0$, $x(0)=\lambda(0)=0$, running the homotopy~\eqref{eq:homotopy} with the parametric active-set strategy of Case~1/Case~2 (together with the dependent-row exchange of Lemma~\ref{lem:dependent-exchange} when needed) reaches, after a finite number $K\le K_{\rm max}$ of steps, where $K_{\rm max}\triangleq\sum_{r=0}^{\min\{m,n\}}\binom mr$, a breakpoint $t_K$ that is either $\ge1$ or marks infeasibility. If $t_K\ge1$, then $x^\star\eqdef x(1)$ is the unique v-GNE of the target problem~\eqref{eq:player}, with associated common multiplier $\lambda^\star\eqdef \lambda(1)$;
otherwise $\{x\in\rr^n:Ax\le b\}=\varnothing$ and~\eqref{eq:player} has no solution.
\label{thm:gPQP}
\end{theorem}

\begin{proof}
We prove by induction that the affine pair~\eqref{eq:affine-path}
associated with $J_j$ solves GNEP$(t)$ on $[t_j,t_{j+1}]$, with
\[
 A_{\bar J_j}x(t)<b_{\bar J_j}(t),\qquad
 \lambda_{J_j}(t)>0,
 \qquad t\in(t_j,t_{j+1}).
\]
For $j=0$, this holds locally to the right of $t_0=0$ by
Proposition~\ref{prop:start}. More generally, suppose these strict
inequalities hold on some interval $(t_j,t_j+\epsilon_j)$. Every
inactive slack $b_i(t)-a_i'x(t)$, $i\in\bar J_j$, and every active
multiplier $\lambda_i(t)$, $i\in J_j$, is affine and strictly positive
immediately to the right of $t_j$. Hence it can first cease to be
positive only at its first root after $t_j$. By~\eqref{eq:events}, $t_{j+1}$ is precisely the
smallest of these roots. Thus the strict inequalities hold throughout
$(t_j,t_{j+1})$ and, by continuity, the corresponding nonstrict
inequalities hold at the endpoints. Proposition~\ref{prop:affine}
therefore shows that the affine pair solves GNEP$(t)$ on
$[t_j,t_{j+1}]$.

By conditions~($a$)--($b$), exactly one ordinary event occurs at a
finite breakpoint. If $t_{j+1}<1$, Lemma~\ref{lem:strict-after-pivot}
supplies the induction hypothesis immediately to the right of
$t_{j+1}$ after an independent-row Case~1 pivot or a Case~2 pivot.
For a dependent entering row, case~($i$) of
Lemma~\ref{lem:dependent-exchange}, condition~($c$), and
Lemma~\ref{lem:strict-after-dependent-exchange} give the same
conclusion. If instead case~($ii$) of
Lemma~\ref{lem:dependent-exchange} applies, then $K(t)=\varnothing$
for all $t>t_{j+1}$; since this branch is invoked only while
$t_{j+1}<1$, in particular
$\{x:Ax\le b\}=K(1)=\varnothing$.

Thus, unless infeasibility is detected, the construction continues
until some $t_K\ge1$ (with $t_K=+\infty$ in Case~0). The affine pair
for $J_{K-1}$ then solves GNEP$(t)$ on $[t_{K-1},t_K]$, so evaluation
at $t=1$ gives $x^\star=x(1)$ and the associated multiplier
$\lambda^\star=\lambda(1)$; uniqueness of $x^\star$ follows from
Proposition~\ref{prop:affine}.

Finally, every $A_{J_j}$ has full row rank, hence $|J_j|\le n$ and
there are at most
\[
 K_{\rm max}
 =\sum_{r=0}^{\min\{m,n\}}\binom mr
\]
possible active sets. Lemma~\ref{lem:no-repeat} shows that none is
revisited, so the procedure terminates after at most $K_{\rm max}$
affine pieces.
\end{proof}

\subsection{Dealing with ties}
\label{sec:relaxing-uniqueness}
So far we have assumed unique minimizers in~\eqref{eq:events}, as a tie would leave Case 1's (resp.\ Case 2's) update undetermined. Exact resolution via lexicographic pivoting (the parametric analogue of Bland's rule, as in~\cite{Bes96}) is possible but costly; instead, we follow the heuristic proposed in qpOASES~\cite[Sec.~3.10]{FKPBD14} 
and set
\[
    \dvec_i=\dvec_{\min}\,\frac{m-i}{m-1}+\dvec_{\max}\,\frac{i-1}{m-1},\qquad 0<\dvec_{\min}<\dvec_{\max}
\]
for $i=1,\ldots,m$, $m>1$. Note that here we reuse only this ramp shape, not qpOASES's warm-start-preserving mechanism, as $x(0)=0,\lambda(0)=0$ solve GNEP$(0)$ for any $\dvec>0$.

The coincidence $\widehat t_{j+1}=\widetilde t_{j+1}<\infty$ excluded by condition (b) of Theorem~\ref{thm:gPQP} is also a very rare event that a suitable choice of $\dvec$ makes unlikely to occur in practice, since $\widehat t_{j+1}$ and $\widetilde t_{j+1}$ vary with $\dvec$ and generically do not cross exactly. Should such a coincidence, or an argmin tie in Case 1 or 2 above, nonetheless be detected numerically, a small jitter of $\dvec$ (e.g., perturbing $\dvec_{\min}$ or $\dvec_{\max}$) can be applied and the algorithm restarted from $t=0$.

\subsection{Linear algebra}
\label{sec:linear-algebra}
The $LU$ factorization of $C$ is computed once, cached, and used whenever computing $x=C^{-1}v$ is required for some $v$ by solving two triangular linear systems. 
When solving linear systems as in~\eqref{eq:two-systems} along the sequence of active sets $J_0,J_1,\ldots$, only the right-hand side and the row added to or removed from $A_J$ change at each pivot. We exploit this structure through the Schur complement of $H_J$ with respect to $C$, so that $H_J$ is never assembled or factored as an $(n+|J|)$-square matrix, as follows. Solving $H_Jx=h$, i.e., $x_c+C^{-1}A_J'x_a=C^{-1}h_c$ and $A_Jx_c=h_a$ with
$x_c,h_c\in\rr^n$ and $x_a,h_a\in\rr^{|J|}$, is equivalent to solving the two linear systems
\begin{equation}
    C y = h_c,\quad S_J x_a = A_J y - h_a
\label{eq:solve_HJ}
\end{equation}
and set $x_c = y - W_J x_a$, 
where $S_J\eqdef A_JC^{-1}A_J'$ is the $|J|\times|J|$ Schur complement of $H_J$ with respect to $C$, which is nonsingular as $H_J$ is invertible by Proposition~\ref{prop:bordered},
and $W_J\eqdef C^{-1}A_J'$ is the $n\times|J|$ matrix of columns $w_i\eqdef C^{-1}a_i'$, $i\in J$
($w_i$ can be computed upfront for every row $a_i$ of $A$, $i=1,\ldots,m$, or
on demand the first time $i$ enters the active set). 

The matrix $S_J$ can be updated recursively when adding row $h$ (Case~1) as
\begin{equation}
    S_{J\cup\{h\}}=\begin{bmatrix}S_J&A_Jw_h\\a_h'W_J&a_h'w_h\end{bmatrix}.
\label{eq:S_j-update}
\end{equation}

In Case~2, dropping row $\ell$ can also be performed efficiently
by appropriately swapping the row and column of $S_J$ corresponding to $\ell$ to the last position, so that the leading $(|J|-1)\times(|J|-1)$ block of $S_J$ corresponds to $S_{J\setminus\{\ell\}}$.

To compute  $(h_{1j},g_{1j})$ and $(h_{2j},g_{2j})$ of~\eqref{eq:two-systems} at every interval via~\eqref{eq:solve_HJ},
we need the LU factorization of $S_J$. In our implementation, we recompute $L_JU_J=S_J$ at every removal of active
constraints, while instead we update $L_JU_J$ recursively at every addition.

\emph{Complexity.} Factoring $C$ costs $O(n^3)$ and computing the required columns $w_i$ costs at most $O(n^2m)$. Each subsequent pivot costs $O(n|J|)$ to extend or truncate $S_J$ from the cached columns of $W$, plus $O(|J|^3)$ to refactor it. 
Periodically recomputing $S_J=A_JW_J$ from scratch and refactoring it, which requires $O(n|J|)+O(|J|^3)$ operations, can be useful to guard against roundoff accumulated by the incremental updates~\eqref{eq:S_j-update}.

\section{Warm starting}
\label{sec:warm-start}
Rather than starting from $x(0)=0$, $\lambda(0)=0$, we may have an existing v-GNE solution $(\bar x,\bar\lambda)$ solving~\eqref{eq:player} when the linear term is $\bar c$ and the right-hand side is $\bar b$, where possibly $c\neq \bar c$ and $b\neq \bar b$. Let $W\subseteq\{1,\ldots,m\}$ be an active set at $\bar x$ such that $A_W$ has full row rank, so that $H_W$ in~\eqref{eq:bordered} is nonsingular by Proposition~\ref{prop:bordered}, and $(\bar x,\bar\lambda)$
satisfies
\begin{equation}
    \begin{aligned}
    & C\bar x+\bar c+A_W'\bar\lambda_W=0,\quad A_W\bar x=\bar b_W\\
    &A_{\bar W}\bar x\leq \bar b_{\bar W},\quad \bar\lambda_W\ge0,\quad \bar\lambda_{\bar W}=0.
    \end{aligned}
\label{eq:KKT-bar}
\end{equation}
Clearly, $W\subseteq\{j\in\{1,\ldots,m\}:\ A_j\bar x =\bar b_j\}$.
Let $\bar v\in\rr^m$, $\bar v_W>0$, $\bar v_{\bar W}=0$, be a perturbation of $\bar\lambda$
such that $\bar\lambda+\bar v>0$, and modify the homotopy~\eqref{eq:homotopy} as follows:
\begin{equation}
    \begin{aligned}
    c(t)\eqdef &\bar c-t(\bar c-c)-(1-t)A_W'\bar v_W\\
    b(t)\eqdef &\bar b + d + tp,\quad p\eqdef b-\bar b-d,
    \end{aligned}
\label{eq:cb(t)warm}
\end{equation}
where $d\geq 0$, $d_W=0$, and $d_{\bar W}>0$.
At $t=0$, $c(0)=\bar c-A_W'\bar v_W$ and $b(0)=\bar b+d$, and $(\bar x,\bar\lambda+\bar v)$ solves GNEP$(0)$ with active set $W$, as
\[
C\bar x+c(0)+A_W'(\bar\lambda_W+\bar v_W)
 \quad=\underbrace{\bigl(C\bar x+\bar c+A_W'\bar\lambda_W\bigr)}_{=0\text{ by~\eqref{eq:KKT-bar}}}-A_W'\bar v_W+A_W'\bar v_W=0
\]
and primal feasibility holds because $A_W\bar x=\bar b_W=b_W(0)$ (using $d_W=0$) and, since $d_{\bar W}>0$, $A_{\bar W}\bar x\leq \bar b_{\bar W}<\bar b_{\bar W}+d_{\bar W}=b_{\bar W}(0)$ holds strictly;
moreover $\bar\lambda+\bar v$ is strictly dual feasible on $W$ (since $\bar v_W>0$).
For $t=1$ we recover the original problem~\eqref{eq:player}. The homotopy path-following algorithm can then be applied to this modified homotopy, starting from $(\bar x,\bar\lambda+\bar v)$ and initial active set $W$.

Such a warm-starting technique can be used to solve a sequence of GNEPs~\eqref{eq:player} in which only the linear term $c$ and/or the right-hand side $b$ change, for example, in a game-theoretic model predictive control setting~\cite{DFB26,HB26}, and to deal with merely monotone games, as we will discuss in Section~\ref{sec:merely-monotone}. 

\section{Equality constraints}
\label{sec:equality}
Equality constraints $Ex=e$, $E\in\rr^{q\times n}$, can be eliminated beforehand to get
a system of the form~\eqref{eq:kkt}, to which the algorithm of Section~\ref{sec:active-set-path-following} applies unchanged.
Without loss of generality, we assume that $E$ has full row rank, $\operatorname{rank}(E)=q$ (cf.~\cite[Remark 5.2]{BT26b}), and compute the full QR factorization
\[
    E'=Q\begin{bmatrix}R\\0\end{bmatrix},\qquad Q=\begin{bmatrix}Q_1&Q_2\end{bmatrix}
\]
with $Q\in\rr^{n\times n}$ orthogonal, $Q_1\in\rr^{n\times q}$, $Q_2\in\rr^{n\times(n-q)}$, and $R\in\rr^{q\times q}$ nonsingular upper triangular, so that $E=R'Q_1'$. Let
$y$ be such that $R'y=e$, and set
\begin{equation}
 N\eqdef Q_2,\qquad x_p\eqdef Q_1y
\label{eq:N-xp}
\end{equation}
where $EN=0$, $N'N=I$, and $x_p$ is a particular solution of $Ex_p=e$. Then,
every element $x\in\{x\in\rr^n:\ Ex=e\}$ is given by $x=x_p+Nz$
for $z=N'(x-x_p)$, $z\in\rr^{n-q}$.

The v-GNE $x^\star\in K$ of the game is characterized by the joint KKT conditions
\begin{equation}
\begin{aligned}
 &Cx^\star+c+A'\lambda^\star+E'\nu^\star=0,&&\lambda^{\star\prime}(Ax^\star-b)=0\\
 &Ex^\star=e,\quad Ax^\star\le b,&& \lambda^\star\geq 0.
\end{aligned}
\label{eq:kkt-eq}
\end{equation}
Substituting $x^\star=x_p+Nz^\star$ and left-multiplying the stationarity condition by $N'$, since $N'E'=(EN)'=0$ we get
\begin{equation}
 \widetilde Cz^\star+\widetilde c+\widetilde A'\lambda^\star=0,\quad 
\widetilde Az^\star\le\widetilde b,
\quad \lambda^{\star\prime}(\widetilde Az^\star-\widetilde b)=0 
\label{eq:kkt-eq-reduced}
\end{equation}
where $\widetilde C\eqdef N'CN$, $\widetilde c\eqdef N'(Cx_p+c)$, $\widetilde A\eqdef AN$, and $\widetilde b\eqdef b-Ax_p$, along with $\lambda^\star\geq 0$. Therefore, we can apply the homotopy path-following algorithm of Section~\ref{sec:active-set-path-following} to the reduced problem~\eqref{eq:kkt-eq-reduced}, which has no equality constraints, to find $z^\star$ and $\lambda^\star$,
and, finally, reconstruct the v-GNE of the original problem as $x^\star=x_p+Nz^\star$. 
Note that the reduced Hessian $\widetilde C$ is still strongly monotone, with the same constant $\mu$ as $C$, since for any $z\in\rr^{n-q}$ we have $z'\widetilde Cz=z'N'CNz=(Nz)'C(Nz)\ge\mu\|Nz\|_2^2=\mu\|z\|_2^2$ because $N$ has orthonormal columns.

If $\nu^\star$ is also needed, it can be recovered by left-multiplying the stationarity condition in~\eqref{eq:kkt-eq} by $E$ to get $\nu^\star=-(EE')^{-1}E(Cx^\star+c+A'\lambda^\star)$ or, equivalently, by solving
the triangular system
$R\nu^\star=-Q_1'\bigl(Cx^\star+c+A'\lambda^\star\bigr)$,
since $E=R'Q_1'$ and, therefore, $(EE')^{-1}E=R^{-1}(R')^{-1}R'Q_1'=R^{-1}Q_1'$.

The elimination approach above requires a one-time $O(nq^2)$ QR factorization and the computation of the reduced matrices needed in~\eqref{eq:kkt-eq-reduced}. 
Similarly to~\cite[Sec.~3.1]{FKPBD14}, which treats equality constraints as box constraints with equal lower and upper bounds, an alternative approach would be to consider such constraints as permanently active and ignore the nonnegativity of their multipliers $\nu$. The drawback would be that maintaining the factorization of $S_J$ would become computationally more expensive, as the active set $J$ would always include the $q$ equality constraints. 

Our implementation of gPAQ follows a third approach: it keeps $E$ in the bordered matrix $\left[\begin{smallmatrix}C&E'\\E&0\end{smallmatrix}\right]$, which is nonsingular by Proposition~\ref{prop:bordered} and is factored once, and uses it in place of $C$ in the Schur-complement scheme of Section~\ref{sec:linear-algebra}. In this way, $S_J$ only involves the active inequality constraints, and the multiplier $\nu$ is obtained directly.
It is easy to show that this implementation inherits all the properties of the elimination approach. The homotopy is run with right-hand side $e(t)\eqdef te$ of the equality constraints in the cold start, and $e(t)\eqdef\bar e+t(e-\bar e)$, $\bar e\eqdef E\bar x$, when warm starting as in~\eqref{eq:cb(t)warm}.

\section{Merely monotone case}
\label{sec:merely-monotone}
Proximal point iterations have been used to deal with convex quadratic programs
in which the Hessian matrix is ill-conditioned or even singular~\cite{FKPBD14,Bem18,ABA22b}.
Similarly, we show in this section that proximal point iterations can be applied to the GNEP~\eqref{eq:player} when $C$ is merely monotone, i.e., $C_s\succeq0$ but Assumption~\ref{ass:strong} does not hold, or $C$ is ill-conditioned due to $\mu\ll1$.

Let us consider the following proximal point iterations in which the cost function of each player in the GNEP~\eqref{eq:player} is modified
as
\begin{equation}
\begin{aligned}
 \min_{x_i}& \frac12 x_i'C_{ii}x_i +x_i'\sum_{j\ne i}C_{ij}x_j+c_i'x_i+\frac{\rho}{2}\|x_i-x_{k,i}\|_2^2\\
\textrm{s.t. }& Ax\le b
\end{aligned}
\label{eq:player-prox}
\end{equation}
where $\rho>0$ is the proximal point parameter acting as a regularization parameter on the GNEP, $x_k$ is the previous iterate, and $x_{k,i}$ denotes its $i$th block. The proximal terms $\frac{\rho}{2}\|x_i-x_{k,i}\|_2^2$ replace $C$ with $C+\rho I$, which makes the regularized game strongly monotone, so that we can apply the parametric active-set method to find the new solution $x_{k+1}$ of~\eqref{eq:player-prox}.

When solving each regularized game~\eqref{eq:player-prox} for $k=0,1,\ldots$ (Hessian $C+\rho I$ and linear term $c-\rho x_k$, producing $x_{k+1}$), we exploit for $k\ge1$ the warm-starting technique of Section~\ref{sec:warm-start} by setting $\bar x=x_k$, $\bar\lambda=\lambda_k$, and $W=W_k$, where $\lambda_k$ and $W_k$ are the common multiplier and the final active set obtained when computing $x_k$, together with $\bar c=c-\rho x_{k-1}$, that is, the linear term against which $x_k$ was computed, and $\bar b=b$. For $k=0$, with $x_0=0$ and $\lambda_0=0$, we apply the cold-started homotopy~\eqref{eq:homotopy} and the trivial canonical start of Proposition~\ref{prop:start} directly to the first regularized problem (Hessian $C+\rho I$, linear term $c-\rho x_0=c$), which is recovered from~\eqref{eq:cb(t)warm} with $W=\varnothing$, $\bar c=0$, and $\bar b=0$.
This allows us to start the homotopy path-following algorithm from a point that, as $k$ grows, is already close to the solution $x_{k+1}$ of~\eqref{eq:player-prox}, potentially reducing the number of required active-set iterations.

\subsection{Convergence analysis}
We analyze the outer iterations~\eqref{eq:player-prox} under the following assumption, which replaces Assumption~\ref{ass:strong}.
\begin{assumption}[Monotonicity]\label{ass:monotone}
$C_s=\frac12(C+C')\succeq 0$, and the set $X^\star\subseteq K$ of v-GNEs of~\eqref{eq:player} is nonempty.
\end{assumption}
Since the diagonal blocks of $C_s$ are the matrices $C_{ii}$, Assumption~\ref{ass:monotone} implies $C_{ii}\succeq0$; accordingly, in this section the standing assumption $C_{ii}\succ0$ made after~\eqref{eq:player} is relaxed to $C_{ii}\succeq0$, while the regularized players' Hessians $C_{ii}+\rho I$ in~\eqref{eq:player-prox} remain positive definite.
Note that $X^\star\ne\varnothing$ holds, for instance, whenever $K$ is nonempty and bounded~\cite{FP03}.
Under Assumption~\ref{ass:monotone}, $X^\star$ need not be a singleton, but it is closed and convex, as it coincides with the solution set of the monotone variational inequality associated with~\eqref{eq:player}~\cite[Sec.~2.3]{FP03}.
For $x\in K$, let $I(x)\eqdef\{j:\ A_jx=b_j\}$ denote the set of active constraints at $x$, and define the set-valued operator
\begin{equation}
 T(x)\eqdef\begin{cases}\{Cx+c+A_{I(x)}'\lambda_{I(x)}:\ \lambda_{I(x)}\ge0\} & \text{if } x\in K\\ \varnothing & \text{otherwise.}\end{cases}
\label{eq:T}
\end{equation}
By~\eqref{eq:kkt}, the v-GNEs of~\eqref{eq:player} are exactly the zeros of $T$, i.e., $X^\star=T^{-1}(0)$.
Moreover, the stacked pseudogradient of the regularized game~\eqref{eq:player-prox} is $Cx+c+\rho(x-x_k)$, so that, by~\eqref{eq:kkt}, $x_{k+1}$ is the v-GNE of~\eqref{eq:player-prox} if and only if
\begin{equation}
 0\in T(x_{k+1})+\rho(x_{k+1}-x_k)
 \quad\Longleftrightarrow\quad
 x_{k+1}=\Bigl(I+\tfrac{1}{\rho}T\Bigr)^{-1}(x_k).
\label{eq:ppa}
\end{equation}
Since $C+\rho I$ satisfies Assumption~\ref{ass:strong} with $\mu=\rho$, $x_{k+1}$ exists, is unique, and is computed exactly, in a finite number of steps, by the parametric active-set method (Theorem~\ref{thm:gPQP}, under its nondegeneracy conditions), provided that $K\neq\varnothing$.
Hence, \eqref{eq:ppa} is the \emph{exact} proximal point algorithm of~\cite{Roc76} applied to $T$ with constant parameter $\gamma_k=1/\rho$ (denoted by $c_k$ in~\cite{Roc76}), and no inexactness criterion is needed.

We can therefore apply the convergence results of~\cite{Roc76}, which we specialize to the LQ-GNEP~\eqref{eq:player}.

\begin{theorem}\label{thm:ppa}
Let Assumption~\ref{ass:monotone} hold, let $\rho>0$, $x_0\in\rr^n$, and let $\{x_k\}$ be generated by~\eqref{eq:player-prox}. Then:
\begin{itemize}
\item[($i$)] $x_k\to x_\infty$ for some v-GNE $x_\infty\in X^\star$ of~\eqref{eq:player};
\item[($ii$)] if there exists $x^\star\in X^\star$ such that $0\in\operatorname{int}T(x^\star)$, then $X^\star=\{x^\star\}$ and $x_k=x^\star$ for all $k$ sufficiently large, i.e., the proximal point iterations converge in a finite number of steps.
\end{itemize}
\end{theorem}
\begin{proof}
The operator $T$ in~\eqref{eq:T} is maximal monotone, being the operator of the monotone variational inequality on $K$ in~\cite[Eq.~(1.3)]{Roc76}; this is the standing assumption under which the results of~\cite{Roc76} used below hold. By~\eqref{eq:ppa}, $\{x_k\}$ is generated by the exact proximal point algorithm $x_{k+1}=P_k(x_k)$, $P_k=(I+\gamma_kT)^{-1}$, with $\gamma_k=1/\rho$ bounded away from zero, which satisfies criteria (A) and (B) of~\cite{Roc76} with $\epsilon_k=\delta_k=0$. Since $X^\star=T^{-1}(0)\neq\varnothing$ by Assumption~\ref{ass:monotone}, $\{x_k\}$ is bounded and~\cite[Thm.~1]{Roc76} applies, giving ($i$) (weak and strong convergence coincide in $\rr^n$). Under the assumption in ($ii$), \cite[Thm.~3]{Roc76} gives $x_k=x^\star$ for all $k$ sufficiently large; moreover, its proof shows that $T^{-1}(w)=\{x^\star\}$ for all $w$ in a neighborhood of the origin, and in particular $X^\star=T^{-1}(0)=\{x^\star\}$.
\end{proof}

The interior condition in Theorem~\ref{thm:ppa}($ii$) can be stated in terms of the KKT conditions~\eqref{eq:kkt}. Let $I^\star\eqdef I(x^\star)$. By~\eqref{eq:T}, $T(x^\star)$ is the translate by $Cx^\star+c$ of the set $\{A_{I^\star}'\lambda_{I^\star}:\ \lambda_{I^\star}\ge0\}$ of nonnegative combinations of the rows of $A_{I^\star}$, whose interior is nonempty if and only if $\operatorname{rank}A_{I^\star}=n$, and is then equal to $\{A_{I^\star}'\lambda_{I^\star}:\ \lambda_{I^\star}>0\}$. Therefore, $0\in\operatorname{int}T(x^\star)$ if and only if
\begin{equation}
 \operatorname{rank}A_{I^\star}=n\quad\text{and}\quad Cx^\star+c+A_{I^\star}'\lambda^\star_{I^\star}=0\ \text{for some}\ \lambda^\star_{I^\star}>0
\label{eq:sharp}
\end{equation}
that is, $x^\star$ is a vertex of $K$ at which it admits a strictly complementary common multiplier. Condition~\eqref{eq:sharp} implies that the v-GNE is unique (cf.~Theorem~\ref{thm:ppa}($ii$)), but does not require $C$ to be strongly monotone. For instance, in the extreme case $C=0$, the v-GNEs of~\eqref{eq:player} are the optimal solutions of the linear program $\min\{c'x:\ Ax\le b\}$, and~\eqref{eq:sharp} holds if and only if this linear program has a unique optimal solution (cf.~\cite[Prop.~8]{Roc76}). On the other hand, \eqref{eq:sharp} can never hold for a v-GNE in the interior of $K$. When~\eqref{eq:sharp} holds, let $\tau>0$ be the radius of the largest ball centered at the origin and contained in $T(x^\star)=\{A_{I^\star}'\lambda_{I^\star}-A_{I^\star}'\lambda^\star_{I^\star}:\ \lambda_{I^\star}\ge0\}$. Since $0\in T(x^\star)+e$ if and only if $-e\in T(x^\star)$, $\tau$ is the largest radius such that $x^\star$ remains the v-GNE for every perturbation $e$ of the linear term $c$ with $\|e\|_2<\tau$. By~\eqref{eq:ppa}, $\rho(x_k-x_{k+1})\in T(x_{k+1})$, so that, as shown in the proof of~\cite[Thm.~3]{Roc76}, $x_{k+1}=x^\star$ as soon as $\rho\|x_{k+1}-x_k\|_2<\tau$. Moreover, by~\cite[Prop.~1(c)]{Roc76} applied to $x_k$ and the fixed point $x^\star$ of $(I+\frac1\rho T)^{-1}$, $\|x_{k+1}-x_k\|_2\le\|x_k-x^\star\|_2\le\|x_0-x^\star\|_2$. Therefore, $x_{k+1}=x^\star$ as soon as $\rho\|x_k-x^\star\|_2<\tau$ and, in particular, $x_1=x^\star$ whenever $\rho<\tau/\|x_0-x^\star\|_2$ (with $x_0\neq x^\star$).

Theorem~\ref{thm:ppa}($ii$) has the following practical consequences. First, when~\eqref{eq:sharp} holds, the overall scheme, i.e., outer proximal iterations~\eqref{eq:player-prox} and inner parametric active-set path-following, terminates after a finite number of active-set pivots with the exact v-GNE, without any stopping tolerance on the outer loop; in particular, $x_{k+1}=x_k$ certifies $x_k\in X^\star$. 
Second, a smaller regularization parameter $\rho$ triggers the termination condition $\rho\|x_k-x^\star\|_2<\tau$ earlier, down to a single outer iteration for $\rho<\tau/\|x_0-x^\star\|_2$, at the price of worse conditioning of $C+\rho I$ in the inner problems.

When~\eqref{eq:sharp} does not hold, finite termination is not guaranteed, but the iterates still converge by Theorem~\ref{thm:ppa}($i$). Moreover, the polyhedral structure of the LQ-GNEP~\eqref{eq:player} guarantees that the convergence is at least linear, as stated in the following theorem.

\begin{theorem}\label{thm:ppa-linear}
Let Assumption~\ref{ass:monotone} hold, let $\rho>0$, $x_0\in\rr^n$, and let $\{x_k\}$ be generated by~\eqref{eq:player-prox}. Then there exist $a\ge0$ and $\delta>0$, depending only on $C,c,A,b$, such that
\begin{equation}
 T^{-1}(w)\subseteq X^\star+a\|w\|_2\mathcal B\qquad\forall w\in\rr^n:\ \|w\|_2\le\delta
\label{eq:upper-lipschitz}
\end{equation}
where $\mathcal B$ is the closed unit ball, and there exists an index $\bar k$ such that
\begin{equation}
 \operatorname{dist}(x_{k+1},X^\star)\le\frac{a\rho}{\sqrt{1+a^2\rho^2}}\operatorname{dist}(x_k,X^\star),\qquad\forall k\ge\bar k.
\label{eq:linear-rate}
\end{equation}
\end{theorem}
\begin{proof}
The graph of $T$ is the union, over all index sets $I\subseteq\{1,\ldots,m\}$, of the polyhedra $\{(x,Cx+c+A_I'\lambda_I):\ A_Ix=b_I,\ A_{\bar I}x\le b_{\bar I},\ \lambda_I\ge0\}$ (each being the image of a polyhedron in $(x,\lambda_I)$ under a linear map). Hence, $T^{-1}$ is a polyhedral multifunction and, by~\cite{Rob81}, it is locally upper Lipschitz at $0$, which is~\eqref{eq:upper-lipschitz}. The rate~\eqref{eq:linear-rate} then follows from~\cite{Luq84}, which extends~\cite[Thm.~2]{Roc76} to non-unique solutions.
\end{proof}

The rate in~\eqref{eq:linear-rate} improves as $\rho$ decreases. Note that the finite-termination case of Theorem~\ref{thm:ppa}($ii$) is a special case of $a=0$, since then $T^{-1}(w)=\{x^\star\}$ for all $w$ in a neighborhood of the origin. More generally, whenever $a=0$, \eqref{eq:linear-rate} gives $x_{k+1}\in X^\star$ for all $k\geq\bar k$, i.e., finite termination also occurs when $X^\star$ is not a singleton.

\section{Strongly monotone nonlinear case}
\label{sec:nl-monotone}
We now consider GNEPs in which player $i$ solves
\begin{equation}
\begin{aligned}
 \min_{x_i}\ & J_i(x_i,x_{-i})\\
 \textrm{s.t. }& Ax\le b
\end{aligned}
\label{eq:player-nl}
\end{equation}
where each $J_i$ is twice continuously differentiable on an open convex set $\Omega\supseteq K$ and convex in $x_i$, and the shared constraints are affine as in~\eqref{eq:player}. The v-GNEs of~\eqref{eq:player-nl} are the solutions of the variational inequality VI$(K,F)$: find $x^\star\in K$ such that
\begin{equation}
 F(x^\star)'(y-x^\star)\ge0,\qquad\forall y\in K
\label{eq:vi-nl}
\end{equation}
where $F(x)\eqdef(\nabla_{x_1}J_1(x),\ldots,\nabla_{x_N}J_N(x))$ is the pseudogradient of the game and $JF(x)$ denotes its (generally asymmetric) Jacobian. We assume that $F$ is strongly monotone on $\Omega$ with modulus $\mu>0$, i.e., $(F(x)-F(y))'(x-y)\ge\mu\|x-y\|_2^2$ for all $x,y\in\Omega$, which implies that $\frac12(JF(x)+JF(x)')\succeq\mu I$ for all $x\in\Omega$, and in particular for all $x\in K$, and that~\eqref{eq:vi-nl} has a unique solution~\cite{FP03} whenever $K\neq\varnothing$. Strong monotonicity on $K$ alone would not suffice for the Jacobian condition when $K$ has empty interior: for instance, $F(x)=(x_1,-x_2)$ is strongly monotone with $\mu=1$ on $K=\{x\in\rr^2:\ x_2=0\}$, but $JF=\diag(1,-1)$. For the LQ-GNEP~\eqref{eq:player}, $F(x)=Cx+c$ and $JF(x)=C$.

To solve~\eqref{eq:vi-nl}, we adopt the globally convergent Newton method of Taji, Fukushima, and Ibaraki~\cite{TFI93}. At the current iterate $x_k\in K$, the pseudogradient is linearized as $F(x)\approx F(x_k)+JF(x_k)(x-x_k)$, and the Newton point $z_k$ is defined as the solution of the linearized variational inequality
\begin{equation}
 \bigl(F(x_k)+JF(x_k)(z_k-x_k)\bigr)'(y-z_k)\ge0,\qquad\forall y\in K.
\label{eq:lin-vi}
\end{equation}
Problem~\eqref{eq:lin-vi} is exactly the v-GNE problem of the LQ-GNEP~\eqref{eq:player} with
\begin{equation}
 C=C_k\eqdef JF(x_k),\qquad c=c_k\eqdef F(x_k)-JF(x_k)x_k
\label{eq:Ck-ck}
\end{equation}
in which player $i$ minimizes the second-order Taylor expansion of $J_i$ around $x_k$, up to terms that do not depend on $x_i$, so that $C_{k,ii}=\nabla^2_{x_ix_i}J_i(x_k)$ and $C_{k,ij}=\nabla^2_{x_ix_j}J_i(x_k)$, $j\neq i$. Since $C_k$ satisfies Assumption~\ref{ass:strong} with the same $\mu$ (which also gives $C_{k,ii}\succeq\mu I$), the Newton point $z_k$ exists, is unique, and is computed exactly in a finite number of steps by gPAQ (Theorem~\ref{thm:gPQP}, under its nondegeneracy conditions).

To globalize the method, \cite{TFI93} uses the regularized gap function of Fukushima~\cite[Eq.~(8)]{TFI93}
\begin{equation}
\begin{aligned}
 f(x)\eqdef&-F(x)'\bigl(H(x)-x\bigr)-\frac{\alpha}{2}\|H(x)-x\|_2^2\\
 H(x)\eqdef&\arg\min_{y\in K}\Bigl\|y-\Bigl(x-\frac1\alpha F(x)\Bigr)\Bigr\|_2
\end{aligned}
\label{eq:gap}
\end{equation}
with $\alpha>0$. The function $f$ is nonnegative on $K$, and $f(x)=0$ with $x\in K$ if and only if $x$ solves~\eqref{eq:vi-nl}~\cite[Prop.~3.1]{TFI93}. Moreover, $f$ is continuously differentiable with gradient $\nabla f(x)=F(x)-\bigl(JF(x)'-\alpha I\bigr)\bigl(H(x)-x\bigr)$~\cite[Prop.~3.2]{TFI93}. Hence, \eqref{eq:vi-nl} is equivalent to minimizing $f$ over $K$, and $f$ is used as a merit function. Evaluating $f(x)$ requires the Euclidean projection $H(x)$ onto the polyhedron $K$, i.e., the solution of a strictly convex QP.

The overall method, which we call gPAN,
is as follows. Given $x_0\in K$, a tolerance $\epsilon\geq0$, and parameters $\beta,\gamma,\sigma\in(0,1)$, at each iteration $k=0,1,\ldots$:
\begin{enumerate}
 \item if $f(x_k)\le\epsilon$, stop;
 \item compute $C_k$, $c_k$ as in~\eqref{eq:Ck-ck}, solve the linearized problem~\eqref{eq:lin-vi}, i.e., the LQ-GNEP~\eqref{eq:player} with data~\eqref{eq:Ck-ck}, with gPAQ to get $z_k$, and set $d_k\eqdef z_k-x_k$;
 \item if $f(x_k+d_k)\le\gamma f(x_k)$, set $t_k=1$; otherwise, set $t_k=\beta^{\ell}$, where $\ell\geq0$ is the smallest integer such that $f(x_k+\beta^{\ell}d_k)\le f(x_k)+\sigma\beta^{\ell}\nabla f(x_k)'d_k$;
 \item set $x_{k+1}=x_k+t_kd_k$.
\end{enumerate}
Note that $x_{k+1}\in K$, since $x_k,z_k\in K$, $K$ is convex, and $t_k\in(0,1]$. As practical safeguards, our implementation projects $x_0$ onto $K$ if needed and caps the number of outer iterations and of backtracking steps per iteration.

The convergence properties of gPAN follow from~\cite{TFI93}. Let $I^\star\eqdef\{j:\ A_jx^\star=b_j\}$ denote the set of active constraints at the v-GNE $x^\star$ of~\eqref{eq:player-nl}, and let $\mathcal F^\star\eqdef\{x\in K:\ A_{I^\star}x=b_{I^\star}\}$ be the minimal face of $K$ containing $x^\star$. We say that \emph{strict complementarity} holds at $x^\star$ if $x\in K$ and $F(x^\star)'(x-x^\star)=0$ imply $x\in\mathcal F^\star$~\cite[Def.~5.1]{TFI93}.

\begin{theorem}\label{thm:tfi}
Let $F$ be continuously differentiable and strongly monotone on an open convex set $\Omega\supseteq K$ with modulus $\mu>0$, let $K\neq\varnothing$, and let $0<\alpha<2\mu$. Then:
\begin{itemize}
\item[($i$)] for any $x_0\in K$, the sequence $\{x_k\}$ generated by gPAN (with $\epsilon=0$) lies in $K$ and, unless it stops at $x^\star$ after finitely many iterations, converges to the unique v-GNE $x^\star$ of~\eqref{eq:player-nl};
\item[($ii$)] if, in addition, $JF$ is Lipschitz continuous on a neighborhood of $x^\star$ and strict complementarity holds at $x^\star$, then there exists $\bar k$ such that $t_k=1$ for all $k\ge\bar k$, and $\{x_k\}$ converges quadratically to $x^\star$.
\end{itemize}
\end{theorem}
\begin{proof}
gPAN is the algorithm of~\cite[Sec.~4]{TFI93} applied to VI$(K,F)$ in~\eqref{eq:vi-nl} with the merit function~\eqref{eq:gap}, in which the linearized problem~\cite[Eq.~(14)]{TFI93} coincides with~\eqref{eq:lin-vi} and is solved exactly by gPAQ. Since $\alpha<2\mu$, ($i$) follows from~\cite[Thm.~4.1]{TFI93}. Since $K$ is polyhedral, ($ii$) follows from~\cite[Thm.~5.1]{TFI93}.
\end{proof}

The condition $\alpha<2\mu$ also guarantees that $d_k$ is a descent direction for $f$ at $x_k$, namely $\nabla f(x_k)'d_k<-(\mu-\frac{\alpha}{2})\|d_k\|_2^2$, whenever $x_k$ is not a solution of~\eqref{eq:vi-nl}~\cite[Lemma~4.1]{TFI93}, so that the backtracking in Step~3 terminates after a finite number of steps. Regarding the strict complementarity condition in Theorem~\ref{thm:tfi}($ii$), let $\lambda^\star$ be a multiplier satisfying the KKT conditions of~\eqref{eq:vi-nl}, namely $F(x^\star)+A'\lambda^\star=0$, $\lambda^\star\ge0$, $\lambda^{\star\prime}(Ax^\star-b)=0$. Then $F(x^\star)'(x-x^\star)=\sum_{j\in I^\star}\lambda^\star_j(b_j-A_jx)$ for all $x\in K$, so that strict complementarity holds at $x^\star$ in particular when $\lambda^\star_j>0$ for all $j\in I^\star$, which is the usual strict complementarity condition.

From an implementation viewpoint, each gPAQ call in Step~2 is warm-started as in Section~\ref{sec:warm-start}, with $\bar x=x_k$ and $W$ being the active set of the previous linearized problem. Indeed, since $\bar c$ and $\bar b$ in~\eqref{eq:cb(t)warm} are arbitrary, the conditions~\eqref{eq:KKT-bar} can be enforced for the current Hessian $C_k$ by setting, for any $\bar\lambda_W\ge0$, $\bar c\eqdef-C_kx_k-A_W'\bar\lambda_W$, $\bar b_W\eqdef A_Wx_k$, and $\bar b_{\bar W}\eqdef b_{\bar W}$. As the iterates converge, the active set stabilizes, so that only a few active-set changes per Newton iteration are needed. The projections $H(x)$ required to evaluate the merit function~\eqref{eq:gap} are computed with the dual active-set QP solver DAQP~\cite{ABA22b}, warm-started from the active set of the previous projection.
Finally, if $JF(x_k)$ is only monotone, the LQ-GNEP~\eqref{eq:lin-vi} can still be solved by the proximal point iterations of Section~\ref{sec:merely-monotone}; however, the global convergence guarantees of~\cite{TFI93} require strong monotonicity of $F$. Extending the method, together with its convergence guarantees, to merely monotone nonlinear GNEPs is a direction for future research. 

\section{Numerical examples}
\label{sec:examples}
All numerical tests are run in Python 3.11 on a MacBook Pro with Apple M5 Max.
The code to reproduce the reported tests is available at \texttt{github.com/bemporad/gpaq\_examples} and can be run via the implementation of \texttt{gPAQ} and \texttt{gPAN} introduced in the \texttt{nashopt} library~\cite{Bem26}.

\subsection{Random LQ-GNEPs}
\label{sec:random-lqgne}
We consider randomly generated $N$-player strongly monotone LQ-GNEPs with $\bar n=2$ variables per player ($x\in\rr^{N\bar n}$), where the pseudogradient $C$ is shifted so that $\lambda_{\rm min}((C+C')/2)=\mu=0.1$. Each instance has $m=3N$ shared inequalities, $\lfloor N/3\rfloor$ of them active at the solution, and $q=\lfloor N/4\rfloor$ shared equalities, with $A,E$ i.i.d.\ standard normal and right-hand sides set for feasibility. We generate 20 instances for each considered number $N$ of agents.

We compare: \texttt{gPAQ}, the homotopy algorithm of Sections~\ref{sec:active-set-path-following}--\ref{sec:equality},
solving each (strongly monotone) instance in a single homotopy pass and serving as the reference solution; \texttt{dr-daqp}, DAQP's~\cite{ABA22b} native Douglas-Rachford AVI solver~\cite{ABB26} applied directly to the variational inequality; and \texttt{OE}, the Operator Extrapolation method~\cite{Lan2022}, likewise applied directly to the variational inequality via DAQP projections warm-started from the previous active set. Both \texttt{dr-daqp} and \texttt{OE} require strong monotonicity, which holds here by construction. 
We also compare with Lemke's method via \texttt{nashopt}'s interface to the \texttt{quantecon} package, solving the following linear complementarity problem (LCP) reformulation of~\eqref{eq:kkt} under the additional lower bound $x\geq \ell$:
\[
    0 \;\leq\;
    \begin{bmatrix} y \\ \lambda \end{bmatrix}
    \;\perp\;
    \underbrace{\begin{bmatrix} C & \bar A' \\ -\bar A & 0 \end{bmatrix}}_{M}
    \begin{bmatrix} y \\ \lambda \end{bmatrix}
    +
    \underbrace{\begin{bmatrix} C\ell + c \\ \bar b - \bar A\ell \end{bmatrix}}_{q}
    =
    \begin{bmatrix} \mu \\ s \end{bmatrix}
    \;\ge\; 0
\]
where $y=x-\ell$, $\bar A = [A'\ E'\ -E']'$, $\bar b = [b'\ e'\ -e']'$, 
$s=\bar b-\bar Ay-\bar A\ell$, and $\mu$ are the Lagrange multipliers of the lower-bound constraints,
and we set $\ell_i=-10^{2}$ for all $i$.

\begin{table}[t]
  \setlength{\tabcolsep}{6pt}
  \renewcommand{\arraystretch}{1.}
  \centering
{%
  \begin{tabular}{rrr|rrrr}
  \hline
  $N$ & $m$ & $q$ & gPAQ & dr-daqp & OE & Lemke \\
  \hline
  2 & 6 & 0 & \textbf{0.002} & 0.004 & 0.075 & 0.143 \\
  5 & 15 & 1 & \textbf{0.002} & 0.008 & 0.144 & 0.032 \\
  10 & 30 & 2 & \textbf{0.006} & 0.024 & 0.262 & 0.086 \\
  20 & 60 & 5 & \textbf{0.029} & 0.101 & 0.603 & 0.483 \\
  50 & 150 & 12 & \textbf{0.424} & 1.303 & 7.592 & 7.463 \\
  100 & 300 & 25 & \textbf{3.335} & 8.722 & 18.128 & 50.747 \\
  150 & 450 & 37 & \textbf{11.633} & 30.173 & 47.806 & 166.748 \\
  200 & 600 & 50 & \textbf{27.726} & 70.752 & 77.441 & 371.786 \\
  \hline
  \end{tabular}
}
\caption{Random strongly monotone LQ-GNEPs: Mean CPU time (ms) over 20  instances per configuration ($N$ players, $\bar n=2$ variables per player, $m$ shared inequalities with $\lfloor N/3\rfloor$ active, $q$ shared equalities). The results of \texttt{Lemke} are averaged over the instances it solved (20, 16, 13, 9, 2, 2, 2, and 3 out of 20, respectively).}
\label{tab:random-lqgne-time}
\end{table}

\begin{table}[t]
  \setlength{\tabcolsep}{7pt}
  \renewcommand{\arraystretch}{1.}
  \centering
{%
  \begin{tabular}{rrr|rrrr}
  \hline
  $N$ & $m$ & $q$ & gPAQ & dr-daqp & OE & Lemke\\
  \hline
  2 & 6 & 0 & \textbf{1} & \textbf{1} & 9 & 7 \\
  5 & 15 & 1 & \textbf{2} & 17 & 21 & 15 \\
  10 & 30 & 2 & \textbf{4} & 30 & 36 & 27 \\
  20 & 60 & 5 & \textbf{7} & 58 & 63 & 52 \\
  50 & 150 & 12 & \textbf{17} & 160 & 165 & 131 \\
  100 & 300 & 25 & \textbf{34} & 347 & 346 & 262 \\
  150 & 450 & 37 & \textbf{51} & 559 & 563 & 389 \\
  200 & 600 & 50 & \textbf{67} & 760 & 758 & 518 \\
  \hline
  \end{tabular}
}
\caption{Random strongly monotone LQ-GNEPs: Mean number of iterations for the same instances as Table~\ref{tab:random-lqgne-time}.}
\label{tab:random-lqgne-iters}
\end{table}

Our proposed \texttt{gPAQ} method is the fastest method at every $N$ (Table~\ref{tab:random-lqgne-time}) and, as shown in Table~\ref{tab:random-lqgne-iters}, performs the smallest number of iterations (active-set changes), tied with \texttt{dr-daqp} at $N=2$. \texttt{OE} is markedly slower in CPU time at small $N$ (about $19\times$ \texttt{dr-daqp}'s time at $N=2$), but the gap narrows sharply with $N$ and nearly closes by $N=200$. Both \texttt{gPAQ} and \texttt{dr-daqp} provide a GNE with a precision on the order of $10^{-14}$ across the whole sweep, \texttt{OE} achieves a precision on the order of $10^{-13}$, \texttt{Lemke} on the order of $10^{-10}$. Lemke's method solves all instances without equality constraints ($N=2$), but often terminates on a secondary ray when equality constraints are present, which enter the LCP above as pairs of opposite inequalities; for $N\geq 50$, it solves only 2 or 3 of the 20 instances.

\subsection{Merely monotone LQ-GNEPs}
To test the proximal-point iterations described in Section~\ref{sec:merely-monotone}, warm-started in Section~\ref{sec:warm-start}, we consider randomly generated $N$-player GNEPs as described in Section~\ref{sec:random-lqgne} except that they are merely monotone. 

\begin{figure}[t]
\centering
\includegraphics[width=\linewidth]{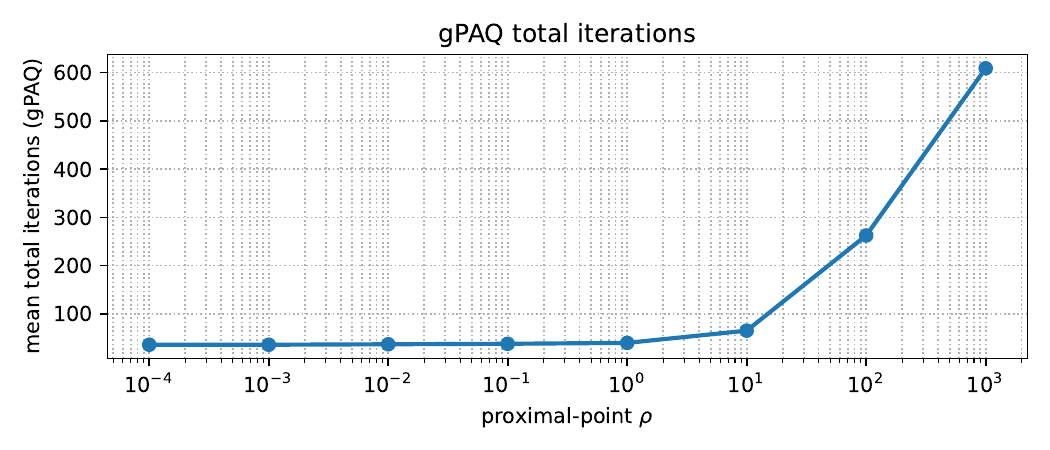}
\caption{Random merely monotone LQ-GNEPs: \texttt{gPAQ}'s mean total iteration count vs the proximal-point regularization $\rho$ ($N=100$, $m=300$, $q=33$, 20 random merely monotone instances per $\rho$).}
\label{fig:merely-rho-sweep}
\end{figure}

First, for $N=100$ agents, we analyze the effect of $\rho$ on \texttt{gPAQ}'s total iteration count
and report the results in Fig.~\ref{fig:merely-rho-sweep} (results are averaged over 20 runs). 
Proximal iterations are stopped when $\rho\|x_{k+1}-x_k\|_2<10^{-10}$, i.e., when the residual $\rho(x_k-x_{k+1})\in T(x_{k+1})$ of the original problem, cf.~\eqref{eq:ppa}, is below $10^{-10}$. The total iteration count is almost constant for $\rho\leq 1$ and grows for larger values of $\rho$, as the contraction factor in~\eqref{eq:linear-rate} increases with $\rho$ and more proximal iterations are needed.

Then, we fix $\rho=0.001$ and compare against \texttt{QP}, the G\"uler-accelerated proximal-point outer loop on the condensed dual QP of~\cite[Section 7.1]{BT26b} ($\rho=10^{-4}$), and \texttt{EG}, Korpelevich's extragradient method~\cite{Kor76} applied directly to the variational inequality, where each iteration uses DAQP~\cite{ABA22b} warm-started from the previous active set. The \texttt{dr-daqp} and \texttt{OE} solvers are excluded, as both require strong monotonicity. Tables~\ref{tab:merely-time}--\ref{tab:merely-iters} report mean CPU time and iteration counts for the three solvers: \texttt{gPAQ} is fastest at every $N$, needing on average only 3--4 proximal iterations regardless of $N$ and a lower overall number of active-set iterations. The proximal iterations do not terminate finitely on these instances but converge linearly, as predicted by Theorem~\ref{thm:ppa-linear}; for $\rho=0.001$ the contraction factor in~\eqref{eq:linear-rate} is so small that each proximal iteration reduces $\|x_k-x^\star\|_2$ by three to five orders of magnitude, so that the stopping criterion is met after very few iterations.

\begin{table}[t]
  \setlength{\tabcolsep}{6pt}
  \renewcommand{\arraystretch}{1.}
  \centering
{%
  \begin{tabular}{rrr|rrr}
  \hline
  $N$ & $m$ & $q$ & gPAQ & QP & EG \\
  \hline
  2 & 6 & 0 & \textbf{0.003} & 0.024 & 7.049 \\
  5 & 15 & 1 & \textbf{0.005} & 0.028 & 2.011 \\
  10 & 30 & 2 & \textbf{0.011} & 0.073 & 1.756 \\
  20 & 60 & 5 & \textbf{0.036} & 0.217 & 2.246 \\
  50 & 150 & 12 & \textbf{0.483} & 3.471 & 10.983 \\
  100 & 300 & 25 & \textbf{3.835} & 28.093 & 40.696 \\
  150 & 450 & 37 & \textbf{12.484} & 95.917 & 98.315 \\
  200 & 600 & 50 & \textbf{30.285} & 230.287 & 191.879 \\
  \hline
  \end{tabular}
}
\caption{Random merely monotone LQ-GNEPs: Mean CPU time (ms) over 20 random merely monotone LQ-GNEP instances per configuration (same $\bar n$, $m$, active constraints, $q$ as in Table~\ref{tab:random-lqgne-time}, $\rho=0.001$).}
\label{tab:merely-time}
\end{table}

\begin{table}[t]
  \setlength{\tabcolsep}{6pt}
  \renewcommand{\arraystretch}{1.}
  \centering
{%
  \begin{tabular}{rrr|rrr}
  \hline
  $N$ & $m$ & $q$ & gPAQ & QP & EG \\
  \hline
  2 & 6 & 0 & \textbf{4 (4)} & 8 & 642 \\
  5 & 15 & 1 & \textbf{4 (4)} & 27 & 179 \\
  10 & 30 & 2 & \textbf{6 (3)} & 51 & 142 \\
  20 & 60 & 5 & \textbf{9 (3)} & 105 & 142 \\
  50 & 150 & 12 & \textbf{19 (3)} & 263 & 136 \\
  100 & 300 & 25 & \textbf{36 (3)} & 552 & 131 \\
  150 & 450 & 37 & \textbf{53 (3)} & 842 & 132 \\
  200 & 600 & 50 & \textbf{69 (3)} & 1131 & 131 \\
  \hline
  \end{tabular}
}
\caption{Random merely monotone LQ-GNEPs: Mean number of iterations for the same instances as Table~\ref{tab:merely-time}: total active-set iterations for \texttt{gPAQ} and \texttt{QP}, extragradient iterations for \texttt{EG}; in parentheses, mean number of outer proximal-point iterations for \texttt{gPAQ}.}
\label{tab:merely-iters}
\end{table}

\subsection{Game-theoretic MPC}
We now turn to structured LQ-GNEPs arising from game-theoretic linear MPC~\cite{BG26,HB26,DFB26}. We consider an MPC problem with $N=2$ agents sharing a linear prediction model, each minimizing a quadratic tracking cost over a prediction horizon $T$, subject to input bounds and \emph{soft} output bounds via slack variables. By eliminating the predicted states from the formulation reported in~\cite[Section 6]{BT26b} and considering output constraints as shared among all agents, we get an LQ-GNEP where each agent decides its own input increments and slack, $z_i=(\delta u_i(0),\ldots,\delta u_i(T-1),\varepsilon_i)$. We choose the MPC weights so that the resulting pseudogradient is strongly monotone. 

For each considered prediction horizon $T$, we draw 20 independent random operating points (initial state, previous input, output reference) and solve the resulting LQ-GNEP with the same three methods as above (\texttt{gPAQ}, \texttt{dr-daqp}, \texttt{OE}); all instances are feasible by construction.

Tables~\ref{tab:mpc-time}--\ref{tab:mpc-iters} report average results for
feasible condensed game-theoretic MPC instances~\cite{DFB26,HB26}.
\texttt{gPAQ} is fastest and requires the fewest iterations in
all cases, while all methods agree to within $2\times10^{-12}$.

\begin{table}[t]
  \setlength{\tabcolsep}{6pt}
  \renewcommand{\arraystretch}{1.}
  \centering
{%
  \begin{tabular}{rrr|rrr}
  \hline
  $T$ & $n$ & $m$ & gPAQ & dr-daqp & OE \\
  \hline
  5 & 22 & 122 & \textbf{0.046} & 0.052 & 4.664 \\
  10 & 42 & 242 & \textbf{0.193} & 0.244 & 28.315 \\
  15 & 62 & 362 & \textbf{0.432} & 0.826 & 43.865 \\
  20 & 82 & 482 & \textbf{0.677} & 1.020 & 43.056 \\
  30 & 122 & 722 & \textbf{1.545} & 2.431 & 130.054 \\
  \hline
  \end{tabular}
}
\caption{Condensed game-theoretic MPC: Mean CPU time (ms) over 20 random operating points per configuration (condensed game-theoretic MPC, $N=2$ agents, $T$ prediction horizon, $n$ decision variables, $m$ shared inequalities).}
\label{tab:mpc-time}
\end{table}

\begin{table}[t]
  \setlength{\tabcolsep}{6pt}
  \renewcommand{\arraystretch}{1.}
  \centering
{%
  \begin{tabular}{rrr|rrr}
  \hline
  $T$ & $n$ & $m$ & gPAQ & dr-daqp & OE \\
  \hline
  5 & 22 & 122 & \textbf{18} & 25 & 449 \\
  10 & 42 & 242 & \textbf{22} & 40 & 1098 \\
  15 & 62 & 362 & \textbf{24} & 77 & 1004 \\
  20 & 82 & 482 & \textbf{21} & 46 & 756 \\
  30 & 122 & 722 & \textbf{22} & 48 & 1208 \\
  \hline
  \end{tabular}
}
\caption{Condensed game-theoretic MPC: Mean number of iterations for the same instances as Table~\ref{tab:mpc-time}.}
\label{tab:mpc-iters}
\end{table}

\subsection{Nonlinear strongly monotone GNEPs}
\label{sec:nl-examples}
We now test the Newton method described in Section~\ref{sec:nl-monotone} on the same random instances used in Section~\ref{sec:random-lqgne}, in which we add to every agent's quadratic cost $J_i$ the jointly convex nonlinear potential term $\Psi(x)=\log\Big(1+\sum_{j=1}^N \exp(a_j'x_j)\Big)$,
where the vectors $a_j$ are drawn i.i.d. from $\mathcal N(0,I)$, and remove equality constraints ($q=0$).
Since $\Psi$ is convex, the pseudogradient Jacobian becomes $C+\nabla^2\Psi(x)\succeq C\succeq\mu I$.
We denote our solution method as \texttt{gPAN}, which linearizes the pseudogradient at each outer Newton iteration and solves the resulting LQ-GNEP with \texttt{gPAQ}. 

We compare \texttt{gPAN} against two general-purpose solvers, denoted by \texttt{trf} (nonlinear least-squares on the KKT residual of the GNEP solved by \texttt{nashopt}) and \texttt{OE}~\cite{Lan2022}, all with stopping tolerance equal to $10^{-10}$. Tables~\ref{tab:nl-time}--\ref{tab:nl-iters} report the mean CPU time and number of iterations, respectively, over 20 instances per configuration (we removed from the reported CPU time the JAX-JIT compile time for each solver). The GNEs found by the three solvers differ on the order of $10^{-8}$ across all instances. Clearly, \texttt{gPAN} is by far the fastest solver.

\begin{table}[t]
  \setlength{\tabcolsep}{6pt}
  \renewcommand{\arraystretch}{1.}
  \centering
{%
  \begin{tabular}{rrr|rrr}
  \hline
  $N$ & $m$ & $m_{\rm act}$ & gPAN & trf & OE \\
  \hline
  2 & 6 & 0 & \textbf{0.574} & 2.969 & 424.832 \\
  3 & 9 & 1 & \textbf{0.469} & 1.970 & 228.885 \\
  5 & 15 & 1 & \textbf{0.557} & 3.595 & 321.374 \\
  10 & 30 & 3 & \textbf{0.725} & 58.735 & 260.457 \\
  \hline
  \end{tabular}
}
\caption{Nonlinear monotone GNEPs: Mean CPU time (ms), excluding each solver's own JIT compile time, over 20 instances per configuration ($N$ players, $\bar n=2$ variables per player, $m$ shared inequalities with $m_{\rm act}$ active at the GNE).}
\label{tab:nl-time}
\end{table}

\begin{table}[t]
  \setlength{\tabcolsep}{6pt}
  \renewcommand{\arraystretch}{1.}
  \centering
{%
  \begin{tabular}{rrr|rrr}
  \hline
  $N$ & $m$ & $m_{\rm act}$ & gPAN & trf & OE \\
  \hline
  2 & 6 & 0 & 2.9 (3.2) & \textbf{9} & 798 \\
  3 & 9 & 1 & 2.5 (3.6) & \textbf{10} & 400 \\
  5 & 15 & 1 & 2.8 (3.8) & \textbf{11} & 495 \\
  10 & 30 & 3 & 2.6 (5.6) & \textbf{11} & 280 \\
  \hline
  \end{tabular}
}
\caption{Nonlinear monotone GNEPs: Mean number of iterations for the same instances as Table~\ref{tab:nl-time}; \texttt{gPAN}'s entries show $a$ ($b$), with $a$ the mean number of outer Newton iterations and $b$ the mean total number of inner \texttt{gPAQ} (homotopy) iterations.}
\label{tab:nl-iters}
\end{table}

\section{Conclusion}
\label{sec:conclusion}
We have shown that Best's parametric active-set method for convex QPs extends to strongly monotone LQ-GNEPs with shared affine constraints, even though the stacked pseudogradient matrix is in general asymmetric. The resulting algorithm, gPAQ, inherits the main features of its QP counterpart, computing the v-GNE exactly, up to round-off errors, after a finite number of active-set changes, or returning a certificate that the shared constraints are infeasible. Unlike other methods, it does not require tuning hyperparameters such as step sizes, and requires very cheap linear algebra. Finally, it can be warm-started from a previously computed solution and active set, which makes it well suited for solving sequences of related games, as in game-theoretic MPC. 

Directions for future research include extending the Newton-type approach and its convergence guarantees to merely monotone nonlinear GNEPs, and exploiting sparsity in the linear algebra to tackle large-scale problems.

\appendix
\section{Appendix}

\subsection{Proof of Lemma~\ref{lem:strict-after-pivot}}
\label{app:strict-after-pivot}
\begin{proof}
Since $A_{J_{j+1}}$ has full row rank,
Proposition~\ref{prop:bordered} ensures that the affine pair
$(\bar x(t),\bar\lambda(t))$ associated with $J_{j+1}$ is well defined
and satisfies~\eqref{eq:active-equalities} for all $t$. It remains to
prove~\eqref{eq:inactive-inequalities}.

Let $W$ be the smaller of $J_j,J_{j+1}$, let $i$ be the index by which
they differ, set $V=W\cup\{i\}$, and write $a=a_i$. Denote by
$(x_W,\lambda_W)$ and $(x_V,\lambda_V)$ the corresponding affine pairs,
and by $\nu(t)$ the multiplier of $i$ on the $V$-branch. Subtracting
their KKT equalities~\eqref{eq:stationarity-t} gives
\[
 H_W
 \begin{bmatrix}
 x_W(t)-x_V(t)\\
 \lambda_W(t)-\lambda_V(t)|_W
 \end{bmatrix}
 =
\begin{bmatrix}a\\0\end{bmatrix}\nu(t).
\]
Thus, if $H_Wz=[a'\ 0']'$, then
$x_V(t)=x_W(t)-\nu(t)z_x$, and using $a'x_V(t)=b_i(t)$ yields
\begin{equation}
 \nu(t)=\frac{a'x_W(t)-b_i(t)}{a'z_x}.
\label{eq:pivot-general}
\end{equation}
By Lemma~\ref{lem:pivot-sign}, $a'z_x>0$: in Case~1 this follows from
the assumed independence of $a_h$ from the rows of $A_{J_j}$, while
in Case~2 it follows from the full row rank of $A_{J_j}$ after
removing $a_\ell$.

\emph{Case 1.}
Here $W=J_j$, $i=h$, $x_W=x$, $x_V=\bar x$, and
$\nu=\bar\lambda_h$. Hence~\eqref{eq:pivot-general} gives
\[
 \bar\lambda_h(t)
 =\frac{a_h'x(t)-b_h(t)}{a_h'z_x}.
\]
By~\eqref{eq:primal-event}, the numerator has root $t_{j+1}$ and
positive slope, so $\bar\lambda_h(t_{j+1})=0$ and
$\bar\lambda_h(t)>0$ for $t>t_{j+1}$. Moreover,
$\widehat t_{j+1}<\widetilde t_{j+1}$ and uniqueness of the minimizer
in~\eqref{eq:primal-event} imply
\[
 \lambda_{J_j}(t_{j+1})>0,
 \qquad
 A_{\bar J_j\setminus\{h\}}x(t_{j+1})
 <b_{\bar J_j\setminus\{h\}}(t_{j+1}).
\]
Since $\bar\lambda_h(t_{j+1})=0$, the two affine branches coincide at
the breakpoint, so $\bar x(t_{j+1})=x(t_{j+1})$ and
$\bar\lambda_{J_j}(t_{j+1})=\lambda_{J_j}(t_{j+1})$.
Proposition~\ref{prop:persistence} therefore preserves these strict
inequalities locally to the right of $t_{j+1}$; together with
$\bar\lambda_h(t)>0$, this proves~\eqref{eq:inactive-inequalities}.

\emph{Case 2.}
Here $W=J_{j+1}$, $i=\ell$, $x_W=\bar x$, $x_V=x$, and
$\nu=\lambda_\ell$. From
$x=x_V=x_W-\lambda_\ell z_x$ and
$a_\ell'x(t)=b_\ell(t)$ we obtain
\[
 a_\ell'\bar x(t)-b_\ell(t)
 =\lambda_\ell(t)a_\ell'z_x.
\]
By~\eqref{eq:dual-event},
$\lambda_\ell(t_{j+1})=0$ and $\lambda_\ell(t)<0$ for
$t>t_{j+1}$; hence
$a_\ell'\bar x(t)<b_\ell(t)$ for $t>t_{j+1}$. Furthermore,
$\widetilde t_{j+1}<\widehat t_{j+1}$ and uniqueness of the minimizer
in~\eqref{eq:dual-event} give
\[
 A_{\bar J_j}x(t_{j+1})<b_{\bar J_j}(t_{j+1}),
 \qquad
 \lambda_{J_{j+1}}(t_{j+1})>0.
\]
Since $\lambda_\ell(t_{j+1})=0$, the two branches coincide at the
breakpoint. Proposition~\ref{prop:persistence} therefore preserves
these strict inequalities locally to the right of $t_{j+1}$; together
with the strict inactivity of $\ell$, this proves
\eqref{eq:inactive-inequalities}.
\end{proof}

\subsection{Proof of Lemma~\ref{lem:dependent-exchange}}
\label{app:dependent-exchange}

\begin{proof}
($i$) Since $\xi_k>0$, the identity $a_h=\sum_{i\in J_j}\xi_i a_i$
gives
\[
 a_k=\xi_k^{-1}
 \left(a_h-\sum_{i\in J_j\setminus\{k\}}\xi_i a_i\right).
\]
Hence replacing $a_k$ by $a_h$ preserves the span and cardinality of
the independent rows of $A_{J_j}$, so $A_{J_{j+1}}$ has full row rank.

Let $\theta^\star=\lambda_k(t_{j+1})/\xi_k$. Using
$a_h=\sum_{i\in J_j}\xi_i a_i$ in the stationarity condition for
$J_j$ gives
\[
 -Cx(t_{j+1})-t_{j+1}c
 =\theta^\star a_h+
 \sum_{i\in J_j\setminus\{k\}}
 \bigl(\lambda_i(t_{j+1})-\theta^\star\xi_i\bigr)a_i
\]
where the coefficient of $a_k$ vanishes by the definition of
$\theta^\star$. By minimality of $\theta^\star$ over $\xi_i>0$,
all the displayed coefficients are nonnegative; the same is immediate
when $\xi_i\le0$. Since $a_h'x(t_{j+1})=b_h(t_{j+1})$ and
$A_{J_j}x(t_{j+1})=b_{J_j}(t_{j+1})$, the pair defined in
\eqref{eq:switch} is therefore a KKT pair for $J_{j+1}$.

It remains to show that the dropped constraint $k$ becomes inactive.
The affine continuation associated with $J_{j+1}$ satisfies
\[
 a_h'\bar x(t)=b_h(t),\qquad
 a_i'\bar x(t)=b_i(t),\quad i\in J_j\setminus\{k\}.
\]
Using the expression for $a_k$ above and
$a_h'x(t)=\sum_{i\in J_j}\xi_i b_i(t)$ along the old $J_j$-path yields
\[
 a_k'\bar x(t)-b_k(t)
 =-\xi_k^{-1}\bigl(a_h'x(t)-b_h(t)\bigr).
\]
The term in parentheses vanishes at $t_{j+1}$ and is strictly positive
for $t>t_{j+1}$ by \eqref{eq:primal-event}. Since $\xi_k>0$,
constraint $k$ is tight at $t_{j+1}$ and strictly inactive thereafter.

\medskip
\noindent
($ii$) Along the $J_j$-path,
\[
 a_h'x(t)=\sum_{i\in J_j}\xi_i b_i(t).
\]
Moreover, by \eqref{eq:primal-event},
$a_h'x(t)-b_h(t)$ vanishes at $t_{j+1}$ and has positive slope.
Hence, for every $t>t_{j+1}$,
\[
 b_h(t)<\sum_{i\in J_j}\xi_i b_i(t).
\]
If $y$ satisfied $a_i'y\leq b_i(t)$ for all $i\in J_j$, then, since
$\xi_i\leq 0$,
\[
 a_h'y
 =\sum_{i\in J_j}\xi_i a_i'y
 \geq\sum_{i\in J_j}\xi_i b_i(t)
 >b_h(t)
\]
so $y$ would violate constraint $h$. Thus $K(t)=\varnothing$ for all
$t>t_{j+1}$, and in particular $K(1)=\varnothing$ whenever
$t_{j+1}<1$.
\end{proof}

\subsection{Proof of Lemma~\ref{lem:strict-after-dependent-exchange}}
\label{app:strict-after-dependent-exchange}
\begin{proof}
By part~($i$) of Lemma~\ref{lem:dependent-exchange},
$(\bar x(t_{j+1}),\bar\lambda(t_{j+1}))$ extends via
\eqref{eq:affine-path} with $J=J_{j+1}$ to an affine pair
$(\bar x(t),\bar\lambda(t))$. Moreover, the proof of
Lemma~\ref{lem:dependent-exchange} gives
\begin{equation}
 a_k'\bar x(t)-b_k(t)
 =-\xi_k^{-1}\bigl(a_h'x(t)-b_h(t)\bigr)
 \qquad \forall t\in\rr.
\label{eq:dependent-k-slack}
\end{equation}

\emph{Strict dual feasibility.}
By~\eqref{eq:switch},
$\bar\lambda_h(t_{j+1})=\theta^\star>0$, since
$\lambda_k(t_{j+1})>0$ and $\xi_k>0$.
For $i\in J_j\setminus\{k\}$ with $\xi_i\le0$,
\[
 \bar\lambda_i(t_{j+1})
 =\lambda_i(t_{j+1})-\theta^\star\xi_i
 \ge \lambda_i(t_{j+1})>0.
\]
If instead $\xi_i>0$, then
\[
 \bar\lambda_i(t_{j+1})
 =\xi_i\left(
 \frac{\lambda_i(t_{j+1})}{\xi_i}-\theta^\star
 \right)>0
\]
where strict positivity follows from the assumed uniqueness of $k$ in
\eqref{eq:switch}. Hence
$\bar\lambda_{J_{j+1}}(t_{j+1})>0$.

\emph{Strict primal feasibility.}
Since
\[
 \bar J_{j+1}=(\bar J_j\setminus\{h\})\cup\{k\}
\]
and $\bar x(t_{j+1})=x(t_{j+1})$, uniqueness of the minimizer in
\eqref{eq:primal-event} gives
\[
 a_i'\bar x(t_{j+1})<b_i(t_{j+1}),
 \qquad i\in\bar J_j\setminus\{h\}.
\]
Moreover, $k\in J_j$ implies
$a_k'\bar x(t_{j+1})=b_k(t_{j+1})$. By
\eqref{eq:dependent-k-slack}, and since
$a_h'x(t)-b_h(t)$ vanishes at $t_{j+1}$ and is strictly positive for
$t>t_{j+1}$, we have
\[
 a_k'\bar x(t)<b_k(t),
 \qquad t>t_{j+1}.
\]

Finally, $\bar\lambda_{J_{j+1}}(t)$ and
$b_i(t)-a_i'\bar x(t)$, $i\in\bar J_j\setminus\{h\}$, are affine and
strictly positive at $t_{j+1}$. Proposition~\ref{prop:persistence}
therefore yields some $\epsilon_{j+1}>0$ such that they remain strictly
positive for
$t\in[t_{j+1},t_{j+1}+\epsilon_{j+1})$.
Together with \eqref{eq:dependent-k-slack}, this gives
$\bar\lambda_{J_{j+1}}(t)>0$,  $A_{\bar J_{j+1}}\bar x(t)<b_{\bar J_{j+1}}(t)$,
for all $t\in(t_{j+1},t_{j+1}+\epsilon_{j+1})$, as claimed.
\end{proof}

\subsection{Proof of Lemma~\ref{lem:no-repeat}}
\label{app:no-repeat}
\begin{proof}
For a fixed full-row-rank set $J_j$, Proposition~\ref{prop:bordered}
implies that the affine pair defined by~\eqref{eq:affine-path} is
uniquely determined by $J_j$. Hence, if $J_j$ were active again at a
later parameter value, it would be represented by the same affine
branch.

At $t_{j+1}$, $J_j$ is left either because an inactive constraint
$h\in\bar J_j$ becomes tight (Case~1) or because an active multiplier
$\ell\in J_j$ becomes zero (Case~2). In Case~1, by
\eqref{eq:primal-event}, the affine function
$a_h'x(t)-b_h(t)$ has root $t_{j+1}$ and positive slope, so
\[
 a_h'x(t)>b_h(t),\qquad t>t_{j+1}.
\]
Thus the old $J_j$-branch is primal infeasible for every later $t$;
this is independent of whether the entering row is linearly
independent or handled through Lemma~\ref{lem:dependent-exchange}.

In Case~2, by~\eqref{eq:dual-event}, $\lambda_\ell(t)$ has root
$t_{j+1}$ and negative slope, so
\[
 \lambda_\ell(t)<0,\qquad t>t_{j+1}
\]
and the old branch is dual infeasible for every later $t$.

Therefore $J_j$ cannot be active again on any interval lying strictly
to the right of $t_{j+1}$.
\end{proof}

\end{document}